\input amstex
\input epsf.sty
\magnification=1200
\documentstyle{amsppt}
\pageheight{9.73 true in}
\vcorrection{-0.8 cm}

\TagsOnLeft 
\CenteredTagsOnSplits
\def\Si{ S \left( \infty \right)}
\def\Go{ G^{\sssize O}}
\def\Ge{ G^{\sssize E}}
\def\Gd{ G^{\sssize D}}
\def\Ga{ \Gamma}
\def\Gao{ \Gamma^{\sssize O}}
\def\Gae{ \Gamma^{\sssize E}}
\def\Gad{ \Gamma^{\sssize D}}
\def\Ke{ K^{\sssize E}}
\def\Ko{ K^{\sssize O}}
\def\Kd{ K^{\sssize D}}
\def\GKd{ \left(\Gd,\Kd \right)}
\def\GKo{ \left(\Go,\Ko \right)}
\def\GKe{ \left(\Ge,\Ke \right)}

\def\KGKd{ \Kd\backslash \Gd/\Kd}
\def\KGK{ K\backslash G/K}
\def\lam{\lambda}
\def\tri{\triangledown}
\def\L{\Lambda}
\def\a{\alpha}
\def\be{\beta}
\def\Mu{\text{\rm M}}
\def\ze{\zeta}
\def\si{\sigma}
\def\dt{\delta}
\def\e{\varepsilon}
\def\phm{\phi^\mu}
\def\tphm{{\tilde \phi}{}^\mu}
\def\Nos{ {\Bbb N}/\si}

\def\mt{\tilde m}
\def\pht{\tilde \phi}
\def\Ht{\tilde H}

\def\ompp{\omega_{p_1,p_2}}
\def\Yab{Y_{a,b}}
\def\YabD{Y_{a,b,D}}
\def\HYab{H(Y_{a,b})}
\def\Hab{H_{a,b}}
\def\Wrs{W_{r,s}}
\def\Ydd{Y_{d_1,d_2}}
\def\FijD{F_{i,j,D}}
\def\bG{ \overline{G}}
\def\bK{ \overline{K}}
\def\bGe{ \overline{G}^{\sssize E}}
\def\bKe{ \overline{K}^{\sssize E}}
\def\bY{ \overline{Y}}
\def\Sig{\Sigma}

\def\gZ{ {\frak Z}}
\def\bZ{ {\Bbb Z}}
\def\bN{ {\Bbb N}}
\def\bZp{ {\Bbb Z}_+}
\def\bZz{ {\Bbb Z}\setminus \{0\}}
\def\bC{ {\Bbb C}}
\def\tr{\operatorname{tr}}
\def\1{\operatorname{1}}
\def\C{\operatorname{C}}

\def\Alt{\operatorname{Alt}}
\def\Sym{\operatorname{Sym}}
\def\supp{\operatorname{supp}}
\def\limind{\operatorname{lim}\, \operatorname{ind}\ }
\def\sgn{\operatorname{sgn}}

\def\Stab{\operatorname{Stab}}
\def\Ind{\operatorname{Ind}}
\def\diag{\operatorname{diag}}

\def\Si{S(\infty)}

\def\ops{\hskip1pt}

\def\wh{\,|\,}
\def\la{\langle}
\def\ra{\rangle}

\topmatter
\title On the representations of the infinite symmetric group
\endtitle
\author Andrei Okounkov 
\endauthor
\abstract
We classify all irreducible admissible 
representations of three Olshanski pairs
connected to the infinite symmetric
group $S(\infty)$. In particular, our methods yield 
two simple proofs of the classical
Thoma's description of the characters of $S(\infty)$.
Also, we discuss  a  certain operation 
called mixture of representations which
provides a uniform construction of all
irreducible admissible representations. 
\endabstract  
\thanks
This is the English version of the author's PhD thesis (1995,
Moscow State University).
\endthanks 
\toc
\widestnumber\subhead{2.2}
\head  0. Introduction \endhead
\subhead 0.1 Tame representations and factor-representations\endsubhead
\subhead 0.2 Olshanski pairs and admissible representations\endsubhead
\subhead 0.3 The statement of the problem and of the results\endsubhead
\head 1. Olshanski semigroups \endhead
\subhead 1.1 Definition and Olshanski's theorem \endsubhead
\subhead 1.2 Parameterization of representations \endsubhead
\subhead 1.3 An example: Thoma multiplicativity \endsubhead
\head 2. Classification of irreducible
admissible representations \endhead
\subhead 2.1 Spectra of the operators $A_i$ in 
spherical representations of the pair $\GKd$ \endsubhead
\subhead 2.2 Spectra of the operators $A_i$ in 
admissible  representations of the pair  $\GKd$\endsubhead
\subhead 2.3 The Thoma theorem \endsubhead
\subhead 2.4 Another proof of Thoma theorem \endsubhead
\subhead 2.5 Classification of the irreducible admissible
representations of the pair $\GKd$ \endsubhead
\subhead 2.6 Description of  $\Ke\backslash\Ge/\Ke$ \endsubhead
\subhead 2.7 Spherical representations of the pair $\GKe$ \endsubhead
\subhead  2.8 Classification of irreducible admissible
representations of the pairs $\GKo$, $\GKe$ \endsubhead
\head 3. Construction of representations \endhead
\subhead 3.1  Mixtures of representations \endsubhead
\subhead 3.2 Mixtures and induction \endsubhead
\subhead 3.3 Elementary representations. \endsubhead
\subhead 3.4 Mixtures in the case of $\GKd$ \endsubhead
\head 4. Concluding remarks \endhead
\endtoc
\address
Department of Mathematics, University of Chicago,
5734 University Avenue, Chicago, IL 60637--1546
\endaddress
\email
okounkov\@math.uchicago.edu
\endemail 
\endtopmatter

\document

\head  0.~Introduction \endhead

\subhead 0.1 Tame representations and factor-representations \endsubhead

Denote by $S(n)$ the group of permutations of the set
$\{1,...,n\}$ and let
$$
\Si=\bigcup_n S(n)
$$
be the union of groups $S(n)$ over all $n$. The group $\Si$ is
one of the simplest examples of a {\it wild} group (see, for example,
the book \cite{23} for the definitions of tame and wild groups,
as well as for other basic notions of infinite-dimensional
representation theory). This means  that the study of all
irreducible representations of $\Si$ does not seem to be a
reasonable problem.

There is probably a unique natural simple topology on the group
$\Si$ induced by the weak (as well as by the strong) operator
topology in the representation by permutations of basis vectors
$$
\Si\,\to\,U(l_2)\,,
\tag 0.1
$$
where $U(l_2)$ denotes the group of unitary operators in the
coordinate Hilbert space $l_2$. Denote by
$$
S_n(\infty)\,, \qquad n=1,2,3,\dots,
$$
the subgroup of $\Si$ which fixes the numbers $1,2,...,n$.
The subgroups $S_n(\infty)$ form a fundamental neighborhood
base of identity in this topology. A unitary representation
$$
\pi: \Si\,\to\,U(H)\,, \tag 0.2
$$
where $U(H)$ is the group of unitary
operators in a Hilbert space $H$, is called {\it tame} if it
is continuous with respect to this
topology on $\Si$ and the weak topology on $U(H)$. All
tame representations  were described
by Lieberman \cite{27} (see also \cite{39}). In
particular, it is known that:
\roster
\item any tame representation is a direct sum of
irreducible ones;
\item irreducible tame representations are labeled by all Young
diagrams $\mu$;
\item the $k$th tensor power of the representation (0.1)
decomposes into irreducible representations  labeled by all Young diagrams
$\mu$ such that $k\ge |\mu| > 0$.
\endroster

These properties make the topological group $\Si$ look similar to
a compact group. One expects the group $\Si$ to be a good combinatorial
model of a big infinite-dimensional group. 
It is clear that the supply of tame
representations is much too small for a unitary dual of 
anything truly infinite-dimensional.

Denote by $\overline{\Si}$ the group of all bijections of the set
of natural numbers. One can introduce a similar topology on this
group and in this topology $\Si$ is a dense subgroup. Tame
representations are precisely those representations
of the group $\Si$ that can be
extended by continuity to the entire group $\overline{\Si}$.

We conclude the discussion of tame representations by one more
technically useful definition of a tame representation.
Denote by
$$
H_n = H^{S_n(\infty)}\,, \quad n=1,2,\dots\,,
$$
the subspace of invariants for the action of the group
$S_n(\infty)$ in the representation (0.2). As shown in
\cite{39}, the representation (0.2) is tame if and only if
$$
H=\overline{\bigcup_n H_n} \,.
$$

Another approach to the representation theory of the group $\Si$
focuses on the study of its factor-representations of finite type
(or, more generally, of semi-infinite type). A beautiful theory
of these representations was developed  by Thoma \cite{49}, and
Vershik and Kerov \cite{6,7,9,10}.

Any finite type factor representation $\pi$
is uniquely determined, up to quasi-equivalence, by its
trace, that is, the restriction of the trace in the factor on the
image of the group. The traces of finite type
factor-representations of the group $G$ are precisely the {\it
characters} of the group $G$. A character of the the group $G$
is, by definition, a function $\phi$ on the
group $G$ which is

\roster
\item central, that is, $\phi(g_1g_2)=\phi(g_2g_1)$ for all
$g_1,g_2\in G$;
\item positive definite, that is, for all $g_1,\dots,g_n\in G$ the
matrix $(\phi(g_i g_j^{-1}))$ is Hermitian and non-negatively
definite;
\item indecomposable, that is, it cannot be represented as a sum of
two linearly independent functions satisfying (1) and (2);
\item normalized by $\phi(e)=1$.
\endroster

If the group $G$ is compact then its characters are precisely
the functions
$$
\phi_\pi(g)=\frac {\tr \pi(g)}{\dim \pi} \,,
$$
where $\pi$ runs over the set of equivalence classes of
irreducible representations of the group $G$. The character
theory of compact groups is a classical chapter of representation
theory. The characters of wild groups were the subject of intense
recent studies, see, for example, \cite{3, 6-11, 13-15, 25, 44, 49-52}.

In the paper \cite{49}, Thoma obtained the following
description of all characters of the group $\Si$. The characters
of the group $\Si$ are labeled by a pair of sequences of real
numbers $\{\a_i\}$, $\{\be_i\}$, $i=1,2,\dots$, such that
$$
\gather
\a_1 \ge \a_2 \ge \a_3 \ge \dots > 0\,,  \quad
\be_1 \ge \be_2 \ge \be_3 \ge \dots > 0\,, \\
\sum \a_i+\sum \be_i\le1\,.
\endgather
$$
The value of the corresponding character on a permutation with a
single cycle of length $k$ is
$$
\sum_i\a_i^k+(-1)^{k-1}\sum_i
\be_i^k \,.
$$
Its value on a permutation with several disjoint cycles
equals the product of the values on each cycle. As usual, it
is assumed that an empty product equals
1. In particular, the character of the regular
representation of the group $\Si$ corresponds to the
sequences $\a_i \equiv 0$, $\be_i \equiv 0$.

The heart of the Thoma's proof is the classification of so
called {\it totally positive} sequences. Recall that a sequence
of real numbers $\{a_i\},i=0,1,2,\dots$ is said to be totally
positive if all minors of the following  infinite Toeplitz matrix
$$
\left[
\matrix
a_0 & a_1 & a_2 & a_3 & \hdots \\
0 & a_0 & a_1 & a_2 & \ddots \\
 & 0 & a_0 & a_1& \ddots \\
 & & 0 & a_0 & \ddots \\
 & & & \ddots & \ddots
\endmatrix
\right]
$$
are non-negative. Thoma obtained the description of all
totally positive sequences using some deep results about entire
functions. He established that a totally positive  sequence
$$
a_0, a_1, a_2, \dots
$$
has generating function of the following form
$$
\sum_i a_i t^i =
e^{\gamma t} \prod_i \frac
{1+ \be_i t} {1-\a_i t}
$$
for some non-negative  $\{\a_i\}$, $\{\be_j\}$, $\gamma$ such
that
$$
\sum \a_i+\sum \be_i < \infty \,.
$$
We shall briefly explain below the connection between totally
positive sequences and representation theory and show how simple
representation theoretical
considerations allow one to simplify Thoma's original
argument significantly and, in particular, to avoid
entire functions entirely.

As a matter of fact, Thoma's description of totally positive sequences was
found earlier in the papers \cite{1, 17}. Totally positive
sequences arise in many problems of analysis (such as
approximation theory or small oscillations), geometry
(convex curves), and probability. An important role is played,
for instance, by the following characteristic property of totally
positive sequences: the convolution of an arbitrary  sequence
with a totally positive
sequence  contains no more changes
of sign than the original  sequence. 
This is a classical result by I.~Schoenberg
\cite{43} which he obtained in the course of his studies of
various generalizations of the Descartes rule. Totally positive
sequences and their continuous analogs were studied by
F.~R.~Gantmakher, M.~G.~Krein, D.~Polya, I.~Schoenberg and
his collaborators, A.~Edrei, S.~Karlin, and others. See, for example,
\cite{1, 16-18, 21} and especially \cite{42} where one can find
references to the most recent applications of the theory of
totally positive sequences. Among the papers discussing the
relations between total positivity and representation theory we
mention \cite{6-11, 12}.

An explicit construction of all corresponding
factor-representations of the group $\Si$ was given later by
A.~M.~Vershik and S.~V.~Kerov \cite{6}. A.~M.~Vershik and
S.~V.~Kerov have also found another proof of Thoma theorem based
on the so called {\it ergodic method}; see \cite{5,7} and also
\cite{57}.

For the group $\Si$, the general ergodic method
specializes to the following procedure. One
starts with a sequence $\lam_n$ of Young diagrams such that
$|\lam_n|=n$. Let
$$
\chi_n (g) = \frac { \tr \pi_{\lam_n} (g)}{\dim \pi_{\lam_n}}
$$
be the normalized character of the symmetric group $S(n)$ corresponding to
the irreducible representation $\pi_{\lam_n}$.

We say that a sequence of characters $\chi_n$ converges, as
$n\to\infty$, to a function $\chi$ on the group $\Si$ if
$$
\chi_n(g)\to\chi(g),\quad n\to \infty
$$
for every element $g\in \Si$. Note that the expression
$\chi_n(g)$ makes sense for all sufficiently large $n$.

It follows from the general approximation theorems \cite{5, 12,
33} that every character $\chi$ of the infinite symmetric group
$\Si$ is a limit of a suitable sequence of characters $\chi_n$ of finite
symmetric groups. A.~M.~Vershik and S.~V.~Kerov proved that
a sequence $\chi_n$ has a limit if and only
if the following limits exist:
$$
\alignat 2
&\lim_{n\to\infty} \frac{ i\text{-th row of } \lam_n}
{n} &&= \a_i,\\
&\lim_{n\to\infty} \frac{ i\text{-th column of } \lam_n}
{n} &&= \be_i \,.
\endalignat
$$
If these limits exist then  the characters $\chi_n$ converge to the
Thoma character with parameters $\{\a_i\}$, $\{\be_i\}$. In other
words, the
parameters $\{\a_i\}$, $\{\be_i\}$ have the meaning of asymptotic lengths
of rows and columns of a Young diagram.

\subhead 0.2 Olshanski pairs and admissible representations
\endsubhead

In \cite{32}, G.~Olshanski initiated the study of a more general
class of representations of the
infinite symmetric group $\Si$. Before giving a definition,
let us consider an
example. Let
$$
\pi: \Si \,\to\, U(M)
$$
be a finite type factor-representation of the group $\Si$. Here
$U(M)$ is the group of unitary operators in a finite factor $M$.
Let the Hilbert space $H$ be the completion of $M$ with
respect to the following Hermitian inner product
$$
(A,B)=\tr AB^* \,.
$$
The group
$$
G=\Si\times\Si\,, \tag 0.3
$$
acts in this space by left and right multiplications and this
representation is irreducible. The identity operator
$$
1\in M \subset H
$$
is the unique vector invariant under  the action of
the diagonal subgroup
$$
K=\diag \Si \subset \Si\times\Si\,, \tag 0.4
$$
and the corresponding matrix element
$$
(\pi(g) 1,1)=\tr \pi(g)
$$
is exactly the trace of the factor representation $\pi$.

Let us  check  that the irreducibility
of the action of the group (0.3) and the existence of a vector invariant under
the action of the subgroup (0.4) implies that the action of the
subgroup (0.4) in the space $H$ is tame \cite{32}.
Indeed,  the subgroups
$$
S(n)\times S(n) \quad\text{and}\quad K_n=\diag S_n(\infty)
$$
commute. Therefore, the subspace
$$
1\in \overline{\bigcup_n H^{K_n}}
$$
is $G$-invariant, and also closed and non-trivial, hence, equal to $H$.

As this example suggests, one should study unitary representations of the group
(0.3) such that their restrictions to the subgroup (0.4) are
tame. Such representations are called  {\it
admissible representations of the pair} (0.4) or simply {\it
representations of the pair} (0.4).

Olshanski's general idea was that, in the
infinite-dimensional situation, it takes two groups to
produce a good representation theory. Namely, one should study the unitary
representations not of a single group $G$, but rather unitary
representations of a pair
$$
K\subset G\,, \tag 0.5
$$
where $K$ is a subgroup of $G$ designated to play the role of a  maximal
compact subgroup of $G$. A unitary representation
$$
G\,\to\,U(H)
$$
of the group $G$ is said to be a {\it representation of the
pair} (0.5) if its restriction on the subgroup $K$ belongs to a
given simple class of representations of the group $K$ (for example, the
class of tame representations). Recall that tame representations
do resemble in many aspects representations of a compact group.
In addition to the pair (0.4), two other pairs closely connected
to the group $\Si$ were considered by G.~Olshanski in \cite{32}
and will be studied in the present paper.

Equivalently, admissible representation can be viewed
as continuous unitary representations of a certain
not locally-compact group $\overline{G}$ containing
$G$ as a dense subgroup, see \cite{32} or the next subsection.

The {\it spherical} representations which, by definition, are
irreducible unitary representations of $G$ with a non-zero $K$-fixed
vector, form a distinguished
 subclass of admissible representations.
It is known (see \cite{32} or below) that in case
case $\dim H^K=1$ and also that the above correspondence
between finite factor representations of $\Si$ and
spherical representations of $(G,K)$ is a bijection.

One of the advantages of the class of
admissible representations is that it is closed under all natural
operations on representations such as restriction to the subgroup
$S_n(\infty)\times S_n(\infty)$,
induction  from such subgroups, taking direct sums and tensor
products.

As it was  shown by Olshanski, every representation of the pair
(0.4) is of von Neumann type I. It follows that the study of general
representations of the pair (0.4) can be essentially reduced to
the study of its irreducible representations. Olshanski
obtained this result  using his so called {\it semigroup method}, which is a very
powerful and beautiful
tool for the study of admissible representations \cite{34,32}.
This method is a far reaching  generalization of the
multiplicativity property of characters of $\Si$ found by Thoma
(see also \cite{15, 19-20, 31}). It can be compared to the use of
Hecke algebras in the $p$-adic representation theory;
the subgroups $K_n$ play the role of the principal congruence
subgroups. The main difference between the $p$-adic groups
and the groups like (0.4) is that $K_n$ are not compact.
Consequently, the definition of the convolution product for
$K_n$-biinvariant functions on $G$ involves a certain
limit transition. On the bright side, in this limit, the
multiplication greatly simplifies and one obtains an actual
semigroup and not just a hypergroup.

As another application of the semigroup machinery
one obtains a way of labeling the irreducible  admissible representations.
All irreducible admissible representation
of the pair (0.4) are indexed by continuous parameters (namely,
parameters used in the Thoma theorem) and some discrete parameters.
More precisely, to every element of the set
$\{\a_i\}\cup\{-\be_j\}\cup\{0\}$ one assigns
two Young diagrams in such a way that all but finitely many diagrams
are empty.  Conversely, any such data correspond, in general, to a representation in a
vector space with an invariant sesquilinear scalar product, not
necessarily positive definite. The classification
problem for irreducible admissible representation
this way can be reformulated as finding all
values of the parameters which correspond to unitary representation.
One can compare this with the problem of
describing the unitary highest weight modules of a Lie
algebra. 

We shall give a brief summary of the semigroup method and the
resulting labeling of representations in Section 1.

In the same paper \cite{32}, G.~Olshanski has constructed a large supply
of irreducible admissible representations. His construction
generalizes, on the
one hand, that of factor-representations of the infinite
symmetric group by A.~M.~Vershik and S.~V.~Kerov, and, on the
other hand, it is an infinite-dimensional generalization of
the classical Hermann Weyl's duality for representations in
traceless tensors. Using this construction and the ergodic
method, G.~Olshanski
obtained two-sided estimates for the set of
parameters corresponding to irreducible admissible
representations.

In the present paper we shall prove, and this is our main result,
that the lower bound  from \cite{32} is actually the correct answer.

\subhead 0.3 The statement of the problem and of the results 
\endsubhead

We study irreducible representations of three similar
Olshanski pairs related to the infinite symmetric group.

By definition, a bijection $g:\bZ \to \bZ$ is called {\it finite} if the set
$$
\{i\in \bZ\wh g(i)\ne i\}
$$
is finite. Define a group $\Go$ as the group of all finite
bijections $\bZ \to \bZ$. Set
$$
\Ge=\{g\in \Go\wh g(0)=0\}
$$
$$
\Gd=\{g\in \Go\wh g(0)=0, g(\bN)=\bN\}
$$
The superscripts $O$, $E$, $D$ mean ``odd'', ``even'' and
``double''. We write simply $G$ if the formulas are valid for all
three groups. Set
$$
G(n)=\{g\in G\wh g(i)=i,|i|>n\}
$$
$$
G_n=\{g\in G\wh g(i)=i,|i|\leq n\}
$$
The group $G$ is the union of an increasing chain of its
subgroups $G(n)$
$$
\{e\}=G(0)\subset G(1)\subset G(2) \subset \dots
$$
$$
\bigcup_n G(n)=G.
$$
The subgroups $G_n$ form a decreasing chain of subgroups,
$$
G_0\supset G_1\supset G_2 \supset \dots,
$$
$$
\bigcap_n G_n=\{e\}.
$$
The subgroups $G(n)$ and $G_n$ commute.

There is a natural involution $i\mapsto -i$ on the set $\bZ$.
Using this involution, one can define a subgroup of the group $G$
which shall play the role of a maximal compact subgroup. Let
$$
K=\{g \in G \wh g(-i)=-g(i)\}.
$$
The pairs $(G,K)$ are Gelfand pairs in the sense of \cite{32}.
Define an increasing and decreasing chains of subgroups
$$
K(n)=K \cap G(n),
$$
$$
K_n=K \cap G_n
$$
in the subgroup $K$.

The group $\Gd$ is isomorphic to a direct product of two copies
of the group $\Si$ (one permutes the positive numbers, the other
--- negative ones). The group $\Kd$ is isomorphic to $\Si$.

The groups $\Go$ and $\Ge$ are isomorphic to $\Si$, and the
subgroups $\Ko$, $\Ke$ are isomorphic to a semi-direct product of
$\Si$ and the group $\bZ_2^\infty$. In this case the Olshanski
pairs differ by the way of embedding of the subgroup into the
group (and the resulting representation theory is also different).

We have $\Gd_n\cong\Gd$ and $\Ge_n\cong\Go_n\cong\Ge$. 
Also $\Gd(n)\cong S(n)\times S(n)$, $\Ge(n)\cong S(2n)$, and 
$\Go(n)\cong S(2n+1)$. 

Let $\pi$ be a unitary representation of the group $G$ in a
Hilbert space $H(\pi)$. Denote by $H(\pi)_n$ the subspace formed by
$K_n$-invariant vectors. This subspace is invariant under the
action of the group $G(n)$. Furthermore, the subspaces $H(\pi)_n$
form an increasing sequence of subspaces
$$
H(\pi)_0 \subset H(\pi)_1 \subset H(\pi)_2 \dots
$$
Hence, their union $\cup_n H(\pi)_n$ is an algebraically
invariant subspace in $H(\pi)$. A representation $\pi$ is called
(see \cite{32}) {\it admissible representation of the pair}
$(G,K)$, if
$$
H(\pi)=\overline{\bigcup_n H(\pi)_n}.
$$
It is known \cite{32} that any admissible representation is of von Neumann
type I.

The notion of admissible representation can also be stated in
topological terms. Namely, as shown in \cite{32}, for every
one of the three pairs $(G,K)$ there exists some topological
group $\bG$ which is not locally compact and contains $G$ as
dense subgroup. The admissible representations of the
pair $(G,K)$ are exactly unitary representations of $G$ admitting
an extension by continuity to the group $\bG$. As a result, the
theory of admissible representations is equivalent to the theory of unitary
representations of  non locally-compact groups $\bG$.
Since every admissible representation generates a von Neumann
algebra of type I, non locally-compact topological groups $\bG$
are actually tame groups.

The construction of the group $\bG$ is as follows. Embed the
group $G$ into the group of all (not necessarily finite)
bijections
$$
g: \bZ \to \bZ.
$$
Consider subgroups
$$
\bK_n=
\left\{g|g(i)=-g(-i),g(i)=i,i>n \right\} \cap G\,,
\quad n=0,1,\dots
$$
and let
$$
\bG=G\cdot\bK  \,.
$$
Finally, define a topology on the group $\bG$ such that a
fundamental neighborhood system of unity is formed by the
subgroups $\bK_n$, $n=0,1,\dots$.

The viewpoint of unitary representations of the groups $\bG$ is
convenient for the constructions of representations, see Section 3.
 Still, the viewpoint of admissible
representations of $(G,K)$ pairs is more convenient for the proof
of classification theorems.

The {\bf main result} of this paper is
a complete description of all
irreducible admissible representations of the three pairs
$(G,K)$ given for the pair $\GKd$ in  Theorem 3, Section 2.5,  
and for the other two pairs in Theorem 5, 
Section 2.8.

The sufficiency of the conditions of given in these theorems
was known before; it follows from Olshanski's
explicit construction of irreducible admissible representations
\cite{32}.  We devote the entire Section 3 to the discussion of
the construction of admissible representations. However, in that
section, our
improvement upon \cite{32} is only a very modest one. It
follows from our classification theorems that the representations
constructed in \cite{32} form actually an open subset of the
admissible dual. Here we introduce the notion of a {\it mixture}
of admissible representations which gives a uniform
construction of all irreducible admissible representation.
This mixture is a kind of an induced representation as
discussed in Section 3.2.

Another known result which we discuss at length in the
present paper is the Thoma's description of characters of
$\Si$. Traditionally, this
result was considered as a very hard one. Our methods
allow to give an new simple proof  (see Section 2.3) and
also to very much simplify Thoma's original proof (see
Section 2.4).

The analog of Thoma theorem for the pair $\GKe$ is
obtained in Theorem 4, Section 2.7.

A brief account of the most important results of this paper
was published in \cite{55}.


\head 1.~Olshanski semigroups \endhead

\subhead 1.1~Definition and Olshanski's theorem \endsubhead

In this section we describe semigroups which shall play a prominent
role in the sequel. These semigroups were introduced by
G.~Olshanski in \cite{32}; they are similar to the 
Brauer semigroups.  A more detailed exposition can be
found in \cite{32, 34}. 

Given two disjoint finite sets $S$ and $S'$, we consider the following set.
An element of $B(S,S')$ is by definition the following data
\roster
\item a partition of the set $S\cup S'$ into pairs together
with a nonnegative real number assigned to each pair, and
\item a finite unordered collection of nonnegative real numbers.
\endroster
Geometrically, this data can be visualized as a compact
$1$-dimensional manifold $M$ such that
\roster
\item
the boundary $\partial M$ of $M$ is $S\cup S'$, and
\item each connected component of $M$ is equipped with a 
non-negative real number (which can be thought of as its 
length).
\endroster 
The connected components with boundaries give a 
partition of $S\cup S'$ into pairs and the lengths
of the cycles of $M$ (by which we mean the loops in $M$,
that is, connected components homeomorphic to $S^1$)
give an array of nonnegative numbers.

It is convenient to position the elements of
$S$ and $S'$ on two horizontal lines, one above the other.  
An example of an 
element of $B(S,S')$ with $|S|=|S'|=5$ is depicted in the following
figure:
$${\epsffile{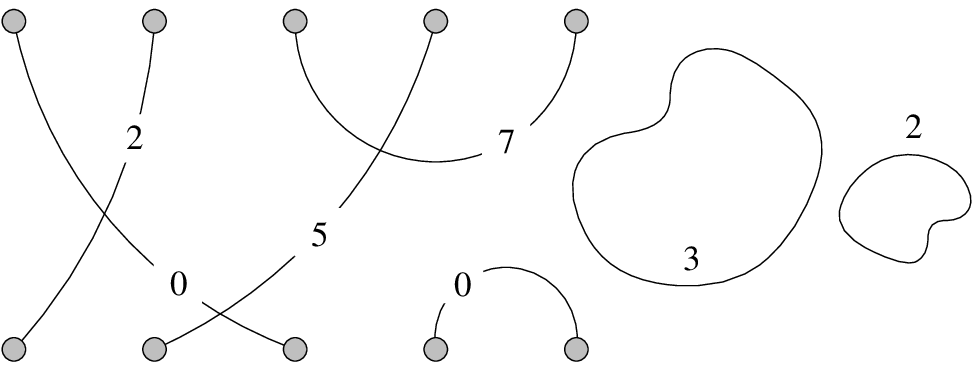}}$$

There is a natural map
$$
B(S,S')\times B(S',S'')\to B(S,S'')
$$
which glues a manifold $M_1\in B(S,S')$ to a manifold
$M_2\in B(S',S'')$ along $S'$. The lengths of the two
glued components, naturally, add up. Note that the 
resulting manifold may have more cycles than $M_1$ and 
$M_2$ combined. This operation makes $B(S,S)$, where
the two $S$'s are considered as two disjoint copies of
the same set, a noncommutative semigroup. For example, if $|S|=5$, then 
$$
\align
&{\epsffile{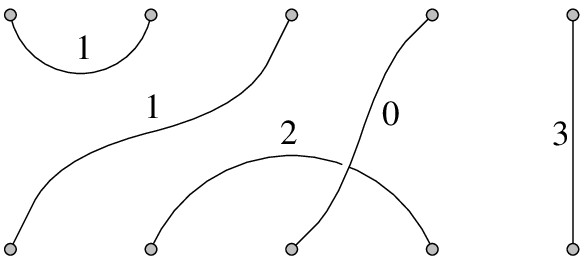}}\\
\vspace{-5\jot}
\intertext{times}  
\vspace{-5\jot}
&{\epsffile{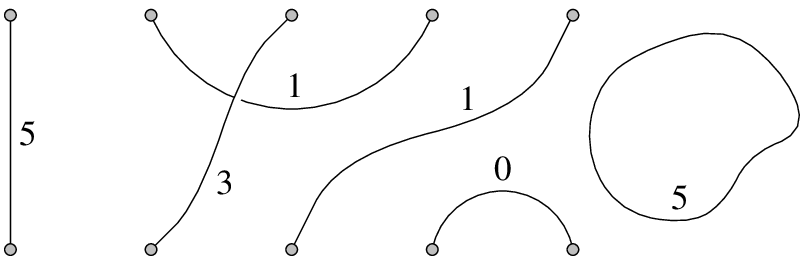}}\\
\vspace{-5\jot}
\intertext{equals}
\vspace{-5\jot}
&{\epsffile{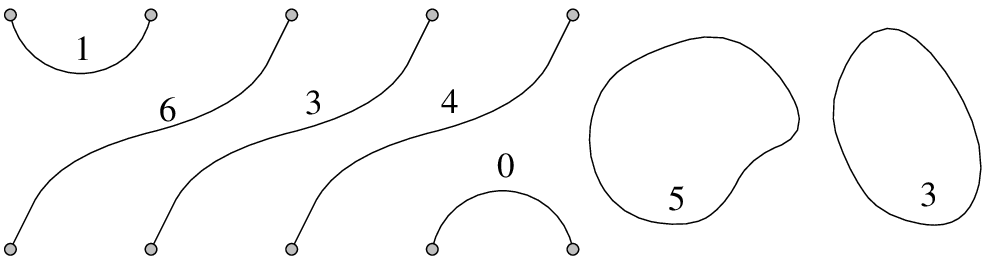}}
\endalign
$$ 
In more technical terms, one can think of an element
of $B(S,S')$ as of a {\it wiring  diagram} with inputs/outputs
indexed by $S\cup S'$ and a nonnegative resistance 
(instead of length) assigned to each wire. 
When one multiplies (that is, connects) two such objects, the resistances
of the connected  wires add up. We shall occasionally
refer to the elements of $B(S,S')$ as (wiring) diagrams.
By analogy with a computer chip, the term {\it chip} was
used in \cite{22,32}.
 
The group $\operatorname{Aut}(S)$ is naturally a subgroup of
$$
\operatorname{Aut}(S)\subset B(S,S)
$$ 
formed by all diagrams with no cycles and all
other components of length zero. Namely, a bijection
$$
g: S \to S
$$
corresponds to such diagram that each element $s\in S$
in the first copy of $S$ is connected to the element
$g^{-1}(s)$ in the second copy of $S$ by a segment of length zero.  
Recall that we consider
the length as just a formal number assigned to a each connected
component; in particular, components of length zero are
still non-trivial. 

Our next goal is to make some sense out of the 
object ``$B(\bZ,\bZ)$''. The above multiplication rule can
fail for infinite diagrams because one can get
infinitely many loops, which is what we want to avoid.

However, we shall need only the semigroup generated by
a certain special set of infinite diagrams; in that semigroup
the multiplication will be indeed well defined. First, we
take all diagrams with no cycles and no components of positive
length. They form a group isomorphic to $\Go$. We add to them
the following diagrams $A_k$, $k\in \bZ$. The diagram 
$A_k$ is defined in the following figure: 
$$
{\epsffile{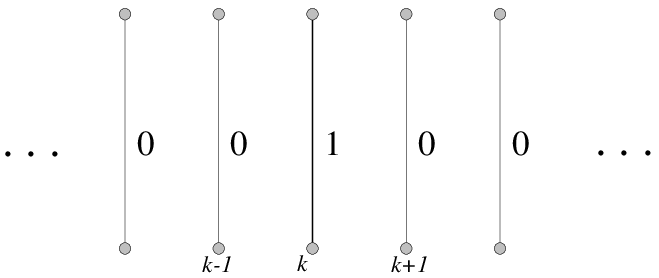}}
$$
All but one segments of $A_k$ have zero length and the only
segment of length $1$ joins $k\in\bZ$ from the first copy of
$\bZ$ with $k$ in the second copy. Clearly,
the semigroup 
$$
\left\langle \Go, A_1 \right\rangle \cong \Si \ltimes 
\bZ_{\ge 0}^{\infty}
$$ 
generated by $\Go$ and $A_1$ contains also all other $A$'s. 
Now, let $C_k$, $k=1,2,3,\dots$, be following pure cycles. The
diagram $C_k$ is, by definition, 
$$
{\epsffile{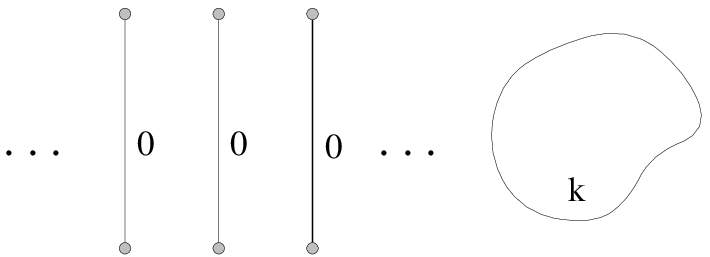}}
$$
In our semigroup, we wish to 
to mod out by the relation 
$$
C_1=1\,.
$$
The purpose of doing this is to make the following diagrams
$P_k$, $k=0,1,2,\dots$, into idempotents. The diagram $P_k$ is 
depicted in the following picture:
$$
{\epsffile{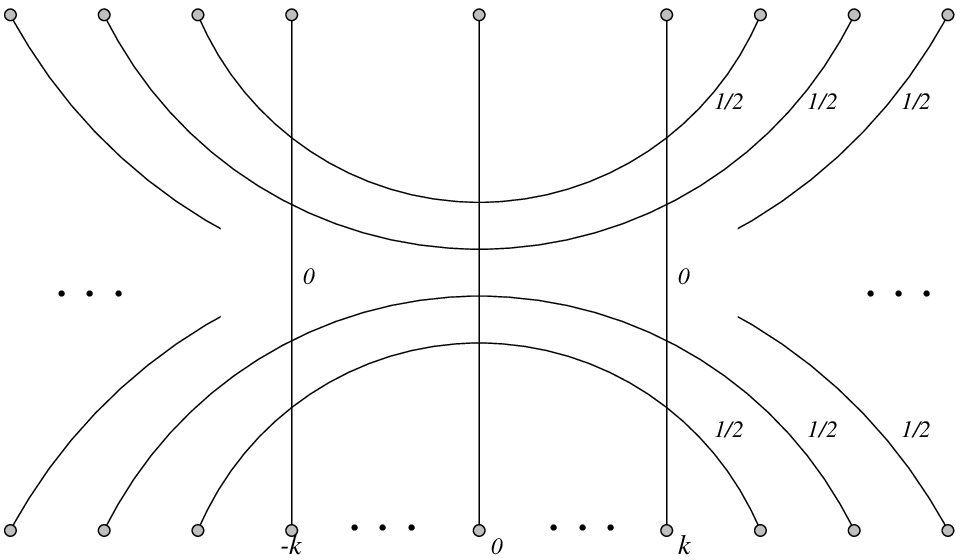}}
$$
It has $2k+1$ vertical segments of length $0$ and all the
remaining components are arcs of length $1/2$. Clearly, $P_k^2=P_k\mod
C_1$.

\definition{Definition}
Denote by $\Gao$ the semigroup generated by the following 
diagrams
$$
\align
\Gao&=\left.\Big\langle \Go, A_1, C_j, P_k \Big\rangle
_{j\ge 1,k\ge 0} \right/ \Big\langle C_1=1 \Big\rangle \\
\intertext{modulo the relation $C_1=1$. Similarly, set}  
\Gae&=\left.\Big\langle \Ge, A_1, C_j, P_k \Big\rangle
_{j\ge 1,k\ge 0} \right/ \Big\langle C_1=1 \Big\rangle\,. \\
\Gad&=\left.\Big\langle \Gd, A_1, C_j, P_k \Big\rangle
_{j\ge 1,k\ge 0} \right/ \Big\langle C_1=1 \Big\rangle\,.
\endalign
$$
\enddefinition

We call these semigroups the {\it Olshanski
semigroups}. They are slightly smaller than the ones defined
by Olshanski 
in \cite{32}, Sections 2.8 and 3.10. Those defined in \cite{32}
are topological
semigroups which contain ours  as dense subsemigroups.

There is a natural involution $*$ in these semigroups, namely,
the reflection in the horizontal axis. This involution
fixes all $A$'s, $C$'s, and $P$'s and takes a permutation
$g$ to the inverse permutation $g^*=g^{-1}$. 

The importance of Olshanski semigroups for the theory
of admissible representations lies in the following
fundamental
\proclaim{Theorem \rm (Olshanski, \cite{32})}
Every admissible representation $\pi$ of a pair $(G,K)$ extends
canonically to a $*$-representation by contractions (that is,
operators of norm $\le 1$) of the
corresponding semigroup $\Ga$ in $H(\pi)$.  In this representation,
the idempotent $P_n$
maps to the orthogonal projection onto $H(\pi)_n$.
\endproclaim

This canonical extension we shall denote with the same letter
$\pi$.

Let us say a few words about how this canonical extension
is constructed. We have to 
specify  the action of the $A$'s and $C$'s. For simplicity,
consider the pair $\GKo$; the other pairs are 
very similar. Consider the transposition $(i,n)\in\Go$
and consider the limit of the operators
$$
\lim_{n\to\infty} \pi((i,n)) \,,
$$
which exists in the weak operator topology. To see this,
 it suffices to check that the limits
$$
\lim_{n\to\infty} (\pi((i,n))\zeta,\eta)
$$
do exist, where the vectors $\zeta$ and $\eta$ belong to the
dense subspace
$$
\bigcup_{m=1}^{\infty} H(\pi)_m \,.
$$
Let $\zeta$ and $\eta$ belong to the subspace $H(\pi)_m$ and
assume that the numbers $n_1$, $n_2$ are chosen big  enough, $n_1,n_2>m$.
Then the permutation
$$
(n_1,n_2)(-n_1,-n_2) \in K_m
$$
belongs to the subgroup $K_m$ and, by the definition of the
subspace $H(\pi)_m$, we have
$$
\multline
\big(\pi((i,n_1))\zeta,\eta\big)=\\
=\big(\pi((n_1,n_2)(-n_1,-n_2)
(i,n_2)(n_1,n_2)(-n_1,-n_2)
)\zeta,\eta\big)= \\
\big(\pi((i,n_2))
\pi((n_1,n_2)(-n_1,-n_2))
\zeta,\pi((n_1,n_2)(-n_1,-n_2))
\eta\big)=\\
=\big(\pi((i,n_2))\zeta,\eta\big)\,.
\endmultline
$$
Therefore, the number
$$
\big(\pi((i,n))\zeta,\eta\big)
$$
does not depend on $n$  provided that $n>m$, hence coincides with
the limit
$$
\big(\pi((i,n))\zeta,\eta\big) =
\lim_{n\to\infty}\big(\pi((i,n))\zeta,\eta\big) \,, \quad n>m \,.
$$
By definition, one sets
$$
\pi(A_i)=\lim_{n\to\infty} \pi((i,n)) \,,
$$
and, similarly,
$$
\pi(C_k)=
\lim_{n_1,\dots,n_k \to\infty}
\pi((n_1,n_2,\dots,n_k)) \,,
$$
where we assume that the numbers $n_1,\dots,n_k$ are pairwise
distinct. 

Sometimes, it is convenient to replace the limits
in the definition of the operators $\pi(A_i)$, $\pi(C_i)$ by the
corresponding Cesaro limits
$$
\pi(A_i)=\lim_{n\to\infty} \frac 1n \sum_{j=1}^n\pi((i,j)) \,,
$$
which exist in the strong operator topology. 
This formula  be interpreted as saying that $A_i$ is the transposition
of $i$ and a ``random'' number $j$. The cycle $C_k$ can be
thought of as a ``random'' cyclic permutation of length $k$. 

The operator $\pi(P_m)$ may be represented in a similar form
$$
\pi(P_m)=\lim_{n\to\infty} \frac 1{n!} \sum_{g\in K_m(n)} \pi(g)\,,
$$
where $K_m(n)=K_m \cap K(n)$.
Indeed, the operator
$$
\frac 1{n!} \sum_{g\in K_m(n)} \pi(g)
$$
is the projection onto the subspace of  $K_m(n)$-invariants. Denote this subspace
by $H(\pi)_{m;n}$. Clearly,
$$
H(\pi)_{m;0} \supset H(\pi)_{m;1} \supset
H(\pi)_{m;2} \supset \dots
$$
and
$$
H(\pi)_m= \bigcap_n H(\pi)_{m;n} \,.
$$
Hence, in the
strong operator topology, the projection onto the subspace $H(\pi)_m$ is the
limit of those onto the subspaces $H(\pi)_{m;n}$.

For the convenience of the future references we list some
useful identities. Given a permutation $\si$,
denote by $[\si]$ the array of
numbers formed by the lengths of non-trivial cycles  of $\si$. One 
easily checks the following equalities: 

\proclaim{Proposition 1}
$$
\alignat3
A_iA_j&=A_jA_i,&&
\tag1.1\\
gA_ig^{-1}&=A_{g(i)},&\qquad&g\in \Go,
\tag1.2\\
A_iP_n&=P_nA_i,&\qquad&|i|\le n,
\tag1.3\\
A_iP_n&=A_{-i}P_n,&\qquad&|i|>n,
\tag1.4\\
P_nA_iP_n&=P_n (i,k) P_n,&\qquad&|i|\le n,\;|k|>n,
\tag1.5\\
P_nA_{i_1}^{k_1}A_{i_2}^{k_2}\dots A_{i_r}^{k_r}P_n
&=P_n\prod_{j=1}^rC_{k_j+1},&\qquad&n<i_j,\;i_m\ne i_l,
\tag1.6\\
P_0\si P_0&=P_0\prod_{k\in[\si]}C_k,&\qquad&\si\in S(\infty)\,.
\tag1.7
\endalignat
$$
\endproclaim

All these identities are straightforward to check. 
As a visual aid for the proof of the last equality we provide
the following figure (where $P_0 (123) P_0= P_0 C_3$ is depicted): 
$$
{\epsffile{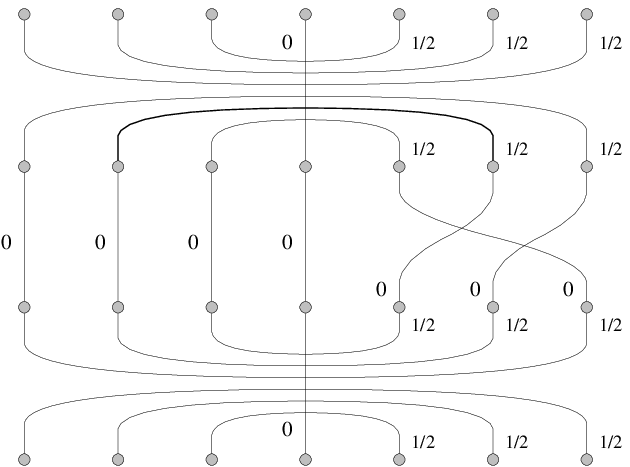}}
$$

\subhead 1.2 Parameterization of representations \endsubhead

Let $\pi$ be an irreducible admissible representation of a pair $(G,K)$ in a
Hilbert space $H(\pi)$. Let $d=d(\pi)$ denote the least integer
$k$, such that $H(\pi)_k \ne 0$; it is called the {\it depth}
of $\pi$. Denote the subspace
$H(\pi)_d$ by $R(\pi)$ and call it the
{\it root} of $\pi$.

Let $\ze$ be a vector in the subspace $R(\pi)$.
The representation $\pi$ is uniquely determined by  any matrix element,
in particular by the one corresponding to $\ze$:
$$
\psi(g)=(\pi(g) \ze,\ze).
$$
The operator $\pi(P_d)$ is the orthogonal projection onto the
subspace $R(\pi)$, hence
$$
(\pi(g) \ze,\ze)=(\pi(P_d g P_d) \ze,\ze).
\tag 1.8
$$
The set
$$
\Ga(d)=P_d \, \Ga \, P_d \subset \Ga 
$$
is a 
subsemigroup  which acts in $R(\pi)$. 
The following proposition is due  to G.~Olshanski:

\proclaim{Proposition 2}
Let $\pi$ be an admissible representation of the group $G$. Choose
$n$ so that $H(\pi)_n\ne 0$. Denote by $\pi_n$ the
representation of the semigroup $\Ga(n)$ in the subspace
$H(\pi)_n$. Then
\roster
\item if the representation $\pi$ is irreducible, then the
representation $\pi_n$ is irreducible;
\item if the representation $\pi_n$ is irreducible and the
subspace $H(\pi)_n$ is cyclic then $\pi$ is also irreducible.
\endroster
\endproclaim

\demo{Proof}
Assume that $\pi$ is irreducible. Let $B$ denote an arbitrary
bounded operator in the subspace $H(\pi)_n$. Denote by $\tilde B$
the operator in the subspace $H(\pi)$ which coincides with $B$ on
the subspace $H(\pi)_n$ and equals zero on its orthogonal
compliment $H(\pi)_n^\perp$. Clearly, this operator is bounded.
Since $\pi$ is irreducible, there exists a sequence $b^{(i)}$ of
elements in the group algebra $\bC[G]$, such that
$$
\pi(b^{(i)}) \to \tilde B, \quad i \to \infty
$$
in the weak operator topology. But this implies that
$$
\pi(P_n b^{(i)} P_n) \to \pi(P_n) \, \tilde B \, \pi(P_n),
\quad i \to \infty \,.
$$
Hence,
$$
\pi_n(P_n b^{(i)} P_n)  \to B\,,
\quad i \to \infty \,.
$$
Therefore, the representation $\pi_n$ is irreducible.

In the opposite direction, let us argue by
contradiction. Let $W$ be a
non-trivial closed invariant subspace. For any $\ze\in H(\pi)_n$,  
 its orthogonal
projections  onto  $W$ and
$W^\perp$ are also $K_n$-invariant vectors. Therefore, at least
one of $\Ga(n)$-invariant subspaces
$$
H(\pi)_n \cap W\ne 0\,, \quad \text{or}\quad H(\pi)_n \cap W^\perp
\ne 0
$$
is non-trivial. Since $\pi_n$ is irreducible, we
conclude  that
$$
H(\pi)_n \subset W\,,  \quad  \text{or}\quad  H(\pi)_n \subset W^\perp\,.
$$
Since the subspace $H(\pi)_n$ is cyclic,
$$
H(\pi) = W\,,\quad \text{or}\quad H(\pi) = W^\perp\,,
$$
which contradicts the non-triviality of $W$.
\qed\enddemo

By virtue of (1.8), we need only  to know the representation
$\pi_d$ of the subgroup $\Ga(d)$ in the subspace $R(\pi)$ in
order to reconstruct the representation $\pi$.

By definition of the number $d$, we have $H(\pi)_{d-1}=0$.
Hence,
$$
\pi(P_{d-1})=0.
$$
Set
$$
\Ga(d)^{\times}=\Ga(d)\setminus \Ga(d) P_{d-1} \Ga(d)\,.
$$
One can check \cite{32} that  $\Ga(d)^{\times}$ is the subsemigroup
of $\Ga(d)$ generated by the elements
$$
\alignat 2
&g P_d\,, \quad &&g\in G(d)\,, \\
&A_i P_d\,, \quad &&|i|\le d \,, \\
&C_k P_d\,, \quad && k\ge 2\,.
\endalignat
$$
It is clear that only the elements of 
this subsemigroups can act in $H(\pi)_d$ by non-zero
operators.

The semigroup $\Ga(d)^{\times}$ is very simple and its
representations can be easily described. In case of $\Gao$,  it is isomorphic to
$$
\Gao(d)^{\times} \cong
(S(2d+1)\ltimes\bZ_+^{2d+1})\times \bZ_+^\infty.
$$
Here the three factors correspond to the three sets of generators of
$\Gao(d)^{\times}$. To simplify notation, set
$$
\Sig(m)=S(m)\ltimes\bZ_+^{m}\,.
$$
All irreducible $*$-representations  of $\Sig(m)$ by contraction
operators are the following. 

Given a point $x\in[-1,1]$ and
a partition $\lam$ of $m$, denote by $V_{\lam,x}$ the 
$\Sig(m)$-module in  which  $\bZ_+^{m}$ acts by  multiplication by $x$ 
and  $S(m)$ acts by the irreducible representation corresponding to $\lam$.
This is an irreducible $*$-representation. More general
representation will be induced from the subsemigroups of the
form
$$
\Sig(\rho)=\prod_i \Sig(\rho_i)\,,
$$
where $\rho$ is a partition of $m$ and the the product is direct. 

By definition, a {\it Young distribution} $\L(x)$ is 
 a function  from
$[-1,1]$ to the set of Young diagrams such that
$\L(x)=\emptyset$  for all but finitely many $x$. We set
$$
\align
|\L|&=\sum_x|\L(x)| \,,\\
\supp\L&=\{x,\L(x)\ne\emptyset\}\,.
\endalign
$$
One can visualize a Young distribution as a collection of Young
diagrams growing out of various points of the interval
$[-1,1]$:
$$
{\epsffile{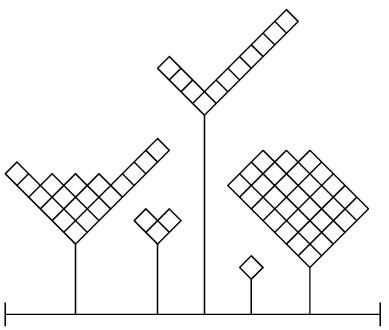}}
$$
By ordering the numbers $|\L(x)|$, $x\in\supp\L$, one
obtains a partition of $|\L|$ which we denote by $\rho(\L)$.
If $|\L|=m$ then
$$
\bigotimes_{x\in\supp\L} V_{\L(x),x}
$$
is an irreducible module over the semigroup
$\Sig(\rho(\L))\subset\Sig(m)$. Denote by $T_\L$ the
representation of $\Sig(m)$ in the following induced
module 
$$
V_\L=\C[\Sig(m)]\otimes_{\C[\Sig(\rho(\L))]}\left( \bigotimes_{x\in\supp\L}
V_{\L(x),x} \right) \,.
$$
This is an irreducible $*$-representation. All irreducible $*$-representations
of $\Sig(m)$ by contractions are precisely the representations in $V_\L$.

The $C$'s are central and $*$-stable, therefore
$$
\pi(C_i)=c_i,
$$
for some  numbers $c_i \in [-1,1]$.

Similarly, we have 
$$
\align
\Gae(d)^{\times} &\cong
(S(2d)\ltimes\bZ_+^{2d})\times \bZ_+^\infty,\\
\Gad(d)^{\times} &\cong
(S(d)\ltimes\bZ_+^{d})\times
(S(d)\ltimes\bZ_+^{d})\times \bZ_+^\infty.
\endalign
$$
The representations $\Gae(d)^{\times}$ are labeled by a Young
distribution $\L$ such that $|\L|=2d$ and by a sequence 
$\{c_i\}\subset [-1,1]$. The representations of
$\Gad(d)^{\times}$ are indexed  by a pair of Young distributions
$\L,\Mu$ such that $|\L|=|\Mu|=d$ and  numbers $c_i$, $i=2,3,\dots$.

The classification problem of admissible irreducible
representations of the group $G$ can now be reformulated as the
problem of description of all parameters $\L$, $\{c_i\}$
(respectively,  $\L$, $\Mu$, $\{c_i\}$) such that an
admissible irreducible representation with these parameters does
exist. Equivalently, one has to describe all parameters for
which the matrix element (1.8) is a positive definite function on
the group $G$.

\subhead 1.3 An example: Thoma multiplicativity \endsubhead

Let $\pi$ be a spherical representation of the pair 
$\GKd$, which is, by definition, an irreducible
unitary representation of $\Gd$ such that $H(\pi)_0\ne 0$. 
By Olshanski's theorem, the
subspace $H(\pi)_0$ carries an irreducible $*$-representation of
the semigroup $\Gad(0)$. This semigroup consists of the
elements of the form
$$
\prod_i C_i^{k_i}P_0.
$$
It is commutative, and has, therefore, only $1$-dimensional 
irreducible representations. Whence,  
$$
\dim H(\pi)_0=1.
$$
In other words, there exists a unique, up to a scalar factor,
$\Kd$-invariant vector $\xi$. We normalize it by the condition
$\|\xi\|=1$ and call the {\it spherical}
vector. Such a vector is determined up to multiplication by a
complex number of absolute value 1. The spherical function
$$
\phi_\pi(g)=(\pi(g)\xi,\xi)
$$
does not depend on the choice of the spherical vector.

Since the
spherical vector $\xi$ is $\Kd$-invariant, the spherical function is constant
on the double cosets  $\KGKd$:
$$
\phi_\pi(k_1 g k_2)\!=\!
\left(\pi(k_1 g k_2)\xi,\xi \right)\!=\!
\left(\pi(g)\pi(k_2)\xi,\pi(k_1^{-1})\xi \right)\!=\!
(\pi(g)\xi,\xi)\!=\!\phi_\pi(g)\,,
$$
where $k_1$ and $k_2$ are two arbitrary permutations in  $\Kd$.

Every double coset  in $\KGKd$ intersects with the subgroup $\Si$
by a conjugacy class in $\Si$. If $\si\in\Si$ and $\si$ has just one
non-trivial cycle of length $k$ then
$$
\phi_\pi(\si)=
(\pi(\si)\xi,\xi)=(\pi(P_0\si P_0)\xi,\xi)\overset(1.7)\to=
(\pi(C_k)\xi,\xi)=(c_k\xi,\xi)=c_k\,.
$$
For arbitrary $\si$, we have
$$
(\si\xi,\xi)=(P_0\si P_0\xi,\xi)\overset(1.7)\to=
\left(\prod_{k\in[\si]}C_k\xi,\xi\right)=\prod_{k\in[\si]}c_k\,.
\tag 1.9
$$
In other words, any spherical function of $\GKd$ is  multiplicative 
in the following sense: its value on any permutation $\si$ equals the product
of its  values on the disjoint cycles of $\si$.
This multiplicativity property was first established by Thoma in \cite{49}.
Olshanski's theorem is, therefore, a generalization of this
Thoma multiplicativity.

\head 2.~Classification of irreducible
admissible representations \endhead

\subhead 2.1 Spectra of the $A_i$'s in 
spherical representations of $\GKd$ \endsubhead

Fix  a spherical representation $\pi$ of the pair $\GKd$ and
write simply $g$ instead of $\pi(g)$.

Let $\mu$ denote the spectral measure for the operator $A_1$ and
the spherical vector $\xi$. Since $\|A_i\|\le 1$,  this measure is
supported at the segment $[-1,1]$.
The numbers $c_k$ are the moments of $\mu$
$$
\int t^k\mu(dt)=(A_1^k\xi,\xi)=(P_0A_1^kP_0\xi,\xi)
\overset(1.6)\to=c_{k+1}\,.
\tag 2.1
$$
Therefore, the spherical function of the pair $\GKd$ is uniquely
determined by the measure $\mu$. We denote the spherical function
corresponding to the measure $\mu$ by $\phm$
$$
\phm(\si)=\prod_{k\in [\si]}{
\int_{[-1,1]} t^{k-1} \, \mu(dt)}.
$$

Let $\lam$ be the measure on $[-1,1]^\infty$ which is the
spectral measure for the operators $A_1,A_2,\dots$ and the vector
$\xi$. The claim of the following lemma is parallel to the
Thoma multiplicativity.

\proclaim{Lemma 1}
$\lam=\mu^{\otimes\infty}$.
\endproclaim

\demo{Proof}
Is suffices to check the identity
$$
\split
\botsmash{\int}t_1^{k_1}\dots t_l^{k_l}d\lam
&=(A_1^{k_1}\dots A_l^{k_l}\xi,\xi)\\
&=(P_0A_1^{k_1}\dots A_l^{k_l}P_0\xi,\xi)
\overset(1.6)\to=\prod_ic_{k_i+1}=\prod_i\int t^{k_i}d\mu
\endsplit
$$
for the integrals of all monomials.
\qed\enddemo

Now we can prove the following

\proclaim{Theorem 1}
The measure $\mu$ is discrete and its atoms can only accumulate
to zero $0\in[-1,1]$.
\endproclaim

\demo{Proof}
Denote by $s$ the transposition $(12)\in \Si$.
Let $E$ be a Borel subset in $[\e,1]$ where $\e>0$. Denote by
$\chi_E$ its characteristic function. We claim that
$$
\e\mu(E)\le\mu(E)^2\,.
$$
To this end, we prove two inequalities
$$
\e\mu(E)\le(s\chi_E(A_1)\ops \xi,\chi_E(A_1)\ops \xi)\le\mu(E)^2.
$$
The expression in the middle is real because
$$
s^{-1}=s\,.
$$
On the one hand, by  (1.3) and (1.5), we have
$$
\split
(s\chi_E(A_1)\ops \xi,\chi_E(A_1)\ops \xi)
&=(s\chi_E(A_1)\ops P_1\xi,\chi_E(A_1)\ops P_1\xi)\\
&\hskip-1pt\overset(1.3)\to=(\chi_E(A_1)\ops P_1sP_1\chi_E(A_1)\ops \xi,\xi)\\
&\hskip-1pt\overset(1.5)\to=(\chi_E(A_1)\ops P_1A_1P_1\chi_E(A_1)\ops \xi,\xi)\\
&\hskip-1pt\overset(1.3)\to=(A_1\chi_E(A_1)\ops \xi,\xi)=\int_Et\,d\mu\ge\e\mu(E)\,.
\endsplit
$$
On the other hand, since $A_1$ and $A_2$ are commuting
projections,
$$
\split
\chi_E(A_1)\ops s\chi_E(A_1)&=\chi_E(A_1)^2s\chi_E(A_1)^2
\overset(1.2)\to=\chi_E(A_1)^2\chi_E(A_2)\ops s\chi_E(A_1)\\
&\hskip-8pt\overset(1.1,2)\to=\chi_E(A_1)\ops \chi_E(A_2)\ops s
\chi_E(A_1)\ops \chi_E(A_2)\,.
\endsplit
$$
Therefore,
$$
\split
(s\chi_E(A_1)\ops \xi,\chi_E(A_1)\ops \xi)
&=(s\chi_E(A_1)\ops \chi_E(A_2)\ops \xi,\chi_E(A_1)\ops \chi_E(A_2)\ops \xi)\\
&\le(\chi_E(A_1)\ops \chi_E(A_2)\ops \xi,\chi_E(A_1)\ops \chi_E(A_2)\ops \xi)=\mu(E)^2,
\endsplit
$$
where the last step relies on Lemma 1.

It follows from the inequality $\e\mu(E)\le\mu(E)^2$ that either
$\mu(E)=0$ or $\mu(E)\ge\e$. An similar estimate holds for
for $E\subset[-1,-\e]$. This implies that the measure $\mu$ is
discrete. Since $\mu$ is a probability measure, there are  no
more than $1/\e$ of its atoms in the interval $[\e,1]$. This
implies the second claim of the theorem.
\qed\enddemo

We denote by $\supp \mu$ the set of atoms of the measure $\mu$.

\subhead 2.2 Spectra of the $A_i$'s in 
admissible  representations of $\GKd$\endsubhead

Now let $\pi$ denote an irreducible admissible representation of
depth $d>0$ corresponding  to some Young distributions $\L$, $\Mu$
and some numbers $\{c_i\}$.

Let $\mu$ be the the spectral measure
$\mu$ of the operator $A_{d+1}$ with respect to some unit vector
$\ze$ in the subspace $R(\pi)$. The measure $\mu$ is independent
of the choice of $\ze$ because the numbers $\{c_i\}$ are the
moments of  $\mu$. By construction, the measure $\mu$
corresponds to a spherical representation of the pair
$$
(\Gd_d, \Kd_d) \cong \GKd\,.
$$
Therefore, $\mu$  is
discrete. Let $\supp \mu$ denote the set of its atoms.

By definition of the distribution $\L$, the spectrum of any
of the operators $A_i$, $i=1,\dots,d$ in the space $R(\pi)$ is
$\supp\L$.

\proclaim{Proposition 3}
$\supp\L\subset\supp \mu\cup\{0\}$,
$\supp\Mu\subset\supp \mu\cup\{0\}$.
\endproclaim

\demo{Proof}
Take $x\in\supp\L\setminus\{0\}$. Let $\ze$ be a vector in the
subspace $R(\pi)$, such that $\|\ze\|=1$ and $A_1\ze=x\ze$.
Denote by $\dt_x$ the function equal to $1$ at the point $x$, and
to $0$ at all other points. Then
$$
\ze=\dt_x(A_1)\ops \ze\,.
$$
Denote by $s$ the permutation $(1,d+1)\in \Si$. Then
$$
\split
0<|x|&=|(A_1\ze,\ze)|\overset(1.5)\to=|(P_dsP_d\ze,\ze)|
=|(s\ze,\ze)|=|(s\dt_x(A_1)\ops \ze,\dt_x(A_1)\ops \ze)|\\
&\hskip-1pt\overset(1.1,2)\to=|(s\dt_x(A_{d+1})\ops \ze,\dt_x(A_{d+1})\ops \ze)|
\le(\dt_x(A_{d+1})\ops \ze,\dt_x(A_{d+1})\ops \ze)=\mu(x),
\endsplit
$$
i.e., $x\in\supp \mu$. The argument for $\Mu$ is analogous.
\qed\enddemo

This proposition was previously proved in \cite{32, Theorem 4.6}
in a more complicated way.

\subhead 2.3 The Thoma theorem \endsubhead

Let $\pi$ be a spherical representation and $\mu$ the corresponding spectral
measure. Denote by $\a_i$, $-\be_i$, $\a_i>0$,
$\be_i>0$ the  non-zero elements in $\supp \mu$. For $x\ne0$ we
set $\nu(x)=\mu(x)/|x|$. By virtue of (2.1),
$$
\multline
c_k=\sum_i\a_i^{k-1}\mu(\a_i)+\sum_i(-\be_i)^{k-1}\mu(-\be_i) \\
=\sum_i\a_i^k\nu(\a_i)+(-1)^{k-1}\sum_i\be_i^k\nu(-\be_i)
\endmultline \tag 2.2
$$
for every $k>1$.

\proclaim{Theorem 2}
The numbers $\nu(\a_i)$, $\nu(-\be_i)$ are positive integers.
\endproclaim

In the proof of this theorem we shall need the following lemma.
Let $\si\in \Si$ be an arbitrary permutation. Denote by $\Nos$
the set of orbits of  $\si$ on the set $\bN$. For an orbit
$p\in \Nos$, denote by $|p|$ its
cardinality.

\proclaim{Lemma 2}
Let $f_i(t)$, $g_i(t)$, $i=1,2,\dots$ be continuous  functions on
$[-1,1]$, all but finitely many identically equal to $1$. Then
$$
\alignat 3
&\text{\rm{a)}}\quad\bigg(\si\prod_{i=1}^\infty f_i(A_i)\ops \xi,\xi\bigg)
=\prod_{p\in \Nos}\int t^{|p|-1}\prod_{j\in p}f_j(t)\,d\mu,
\\
&\text{\rm{b)}}\quad\bigg(\si\prod_{i=1}^\infty f_i(A_i)\ops \xi,
\prod_{i=1}^\infty g_i(A_i)\ops \xi\bigg)=\prod_{p\in \Nos}\int
t^{|p|-1}\prod_{j\in p}f_j(t)\ops \overline
{g_j(t)}\,d\mu\,,
\endalignat
$$
for $\si\in S(\infty)$.
\endproclaim

\demo{Proof}
a) One can assume that $f_i(t)=t^{k_i}$, $i=1,2,\dots$, and that
$k_i=0$ for $i\gg1$. For $p\in \Nos$ set
$\Sigma(p)=\sum_{j\in p}k_j$. Then $\prod_{j\in p}f_j(t)=t^{\Sigma(p)}$.
The following identity
$$
P_0\si\prod A_i^{k_i}P_0=P_0\prod_{p\in \Nos}C_{|p|+\Sigma(p)}
$$
generalizes the identities (1.6) and (1.7) and is proved
similarly. It follows from this identity that
$$
\split
\bigg(\si\prod_{i=1}^\infty f_i(A_i)\ops \xi,\xi\bigg)
&=\Big(\si\prod A_i^{k_i}\xi,\xi\Big)=\prod_{p\in \Nos}c_{|p|+\Sigma(p)}\\
&=\prod_{p\in \Nos}\int t^{|p|+\Sigma(p)-1}d\mu
=\prod_{p\in \Nos}\int t^{|p|-1}\prod_{j\in p}t^{k_j}d\mu\,.
\endsplit
$$
Part b) follows from a), equation (1.2), and from the obvious
equality
$$
\prod_{j\in p}f_j(t)=\prod_{j\in p}f_{\si(j)}(t)\,. \qed 
$$
\enddemo

\noindent{\bf Remark.}
The above lemma holds  for a larger class of
functions, for example, for functions which are pointwise limits of
uniformly bounded sequences of continuous functions. This follows
from the functional calculus of operators, cf. \cite{41, v.1,
Theorem VII.2(d)}. The function $\dt_x(t)$,
$$
\dt_x(t)=\cases
1\,,  &t=x\,, \\
0\,, &t\ne x\,,
\endcases
$$
belongs to this class. Alternatively, by Theorem 1, we
can take  instead of $\dt_x$ a continuous function
equal to $1$ at the point $x$, and to zero at other points in
$\supp \mu$.

\demo{Proof of Theorem 2}
Fix some $\a=\a_i$ and set $\nu=\nu(\a_i)$. Consider the vector
$$
\smash{\ze^{(m)}=\prod_{i=1}^m\dt_\a(A_i)\ops \xi\,.}
$$
By Lemma 2,
$$
(\si\ze^{(m)},\ze^{(m)})=\prod_{p\in \Nos}\a^{|p|-1}\mu(\a)
=\a^{m-\ell(\si)}\mu(\a)^{\ell(\si)}=\a^m\nu^{\ell(\si)}
$$
for $\si\in S(m)$, where $\ell(\si)$ denotes the number of cycles of
the permutation $\si\in S(m)$. Effectively, what the consideration
of the vectors $\ze^{(m)}$, $m=1,2,\dots$, allows us is to single out
just one point $\a$ from the set $\supp \mu$. 

Let $\Alt\ops (m)$ be the operator of antisymmetrization over the group
$S(m)$. Then
$$
\split
(\Alt\ops (m)\ops \ze^{(m)},\ze^{(m)})
=\frac1{m!}\sum_{\si\in S(m)}\sgn\ops (\si)(\si\ze^{(m)},\ze^{(m)})=\\
=\frac{\a^m}{m!}\sum_{\si\in S(m)}\sgn\ops (\si)\ops \nu^{\ell(\si)}
=\frac{\a^m}{m!}\sum_{\si\in S(m)}(-1)^{m+\ell(\si)}\nu^{\ell(\si)}=\\
=\frac{\a^m}{m!}\,\nu(\nu-1)\dots(\nu-m+1),
\endsplit \tag 2.3
$$
where we have applied the equality $\sgn\si=(-1)^{m+\ell(\si)}$ and
the well known identity
$$
\sum_{\si\in S(m)}x^{\ell(\si)}=x(x+1)\dots(x+m-1)
\tag 2.4
$$
(cf. \cite{45, Proposition 1.3.4}). Since $\Alt\ops (m)$ is a
projection, the last product in (2.3) should be nonnegative for
all $m$ which is only possible if it is terminating, i.e., if
$\nu=\nu(\a)\in\Bbb N$.

If we replace the representation $T$ by its tensor product with  
the representation $\sgn\otimes\sgn$, then
the measure $\mu(x)$ is replaced by the
measure $\mu(-x)$. It follows that $\nu(-\be_i)\in\Bbb N$.
\qed\enddemo

We call  to a discrete probability measure $\mu$ on
$[-1,1]$ satisfying
$$
\frac{\mu(x)}{|x|}\in \bZ_+\,, \quad x\ne 0\,,
$$
a {\it Thoma measure}.

The theorem just proved provides a necessary condition for the
existence of representations. It follows from the explicit
construction of representations \cite{6, 32} (see also Chapter 3),
or else from a direct verification of the positive definiteness
\cite{49}, that this condition is
sufficient as well. Therefore the description of spherical functions
for the pair $\GKd$ is established. In order to state it in the
classical form, it is convenient to treat the set $\supp \mu$ as a
multiset in which every element $x$ is
repeated $\mu(x)/|x|$ times.

\proclaim{Theorem {\rm (Thoma,  \cite{49})}}
The characters of the group $S(\infty)$ are precisely the
functions of the form
$$
\botsmash{
\phi(\si)=\prod_{k\in[\si]}\bigg(\sum_i\a_i^k+(-1)^{k-1}\sum_i
\be_i^k\bigg),}
$$
where
$$
\gather
\a_1 \ge \a_2 \ge \a_3 \ge \dots > 0  \quad
\be_1 \ge \be_2 \ge \be_3 \ge \dots > 0 \\
\sum \a_i+\sum \be_i\le1\,.
\endgather
$$
\endproclaim

\subhead 2.4 Another proof of Thoma theorem \endsubhead

In this Section we show how the presentation
$$
c_k=\int t^{k-1}d\mu
$$
simplifies the original proof of the Thoma theorem. Throughout the
section, except for the very last punch-line, we closely
follow Thoma's original exposition \cite{49}.

Let $\phi$ be a character of  $\Si$. Consider the
restriction of $\phi$ to a finite symmetric group
$S(n)$. The characters of finite symmetric group $S(n)$ are
labeled by Young diagrams with $n$ boxes. Let $\chi^\lam$ be the
(non normalized) character corresponding to an irreducible
representation $\lam$. Since the function $\phi$ is positive
definite, its restriction $\phi|_{S(n)}$ to the group $S(n)$
is a non-negative linear combination of the functions $\chi^\lam$
$$
\phi|_{S(n)}=\sum_{\lam,|\lam|=n} {m(\lam) \, \chi^\lam},
\quad m(\lam) \ge 0 \,.
$$
We call the numbers $m(\lam)$ the {\it Fourier coefficients} of the
function $\phi$. One can compute  them using the orthonormality
of characters with respect to the Hermitian inner product
$$
\la f_1, f_2 \ra_{S(n)}=
\frac 1{n!} \sum_{g\in S(n)} f_1(g)
\overline{f_2(g)}\,.
$$
Evidently, the numbers $m(\lam)$ have to satisfy some coherency
conditions. Namely, consider the representation of the group
$S(n+1)$ determined by a Young diagram $\L$, $|\L|=n+1$.
According to the Young branching rule 
$$
\chi^\L|_{S(n)}=\sum_{\lam,\L\searrow\lam}
\chi^\lam \,,
\tag 2.5
$$
where the notation $\L\searrow\lam$ means that the diagram $\lam$
is obtained from the diagram $\L$ by removing a box. By virtue of
(2.5), the numbers $m(\lam)$ have to  satisfy the conditions
$$
m(\lam)=\sum_{\L,\L\searrow\lam} m(\L)\,.
$$
Conversely, every collection of non-negative numbers
$m(\lam)$ satisfying the above coherence condition determines some
positive definite function on the group $\Si$. The function is
normalized if and only if
$$
m(\emptyset)=1\,.
$$
By Proposition 2 (see also Section 1.3) 
this function is indecomposable if and only if it
is multiplicative in the cycles a of permutation. This
multiplicativity imposes severe restrictions on
the numbers $m(\lam)$. Our next goal is to
obtain a precise form of these restrictions.

Recall the definition of the external product of
characters of symmetric groups. Let $\gZ(S(n))$ denote
the linear space of central functions on the groups $S(n)$
and suppose $f_1\in \gZ(S(n_1))$ and $f_2\in \gZ(S(n_2))$ are
some central functions. Consider the following function
on $S(n_1)\times S(n_2)$ 
$$
(f_1 \otimes f_2)(g_1,g_2)=f_1(g_1) f_2(g_2)\,.
$$
The external product of $f_1$ and $f_2$ is, by definition, 
$$
f_1 \circ f_2 = \Ind_{\,S(n_1)\times S(n_2)}
^{\,S(n_1+n_2)} f_1 \otimes f_2 \quad  \in \gZ(S(n_1+n_2)) \,.
$$
Denote by $\eta^k$ the trivial character of $S(k)$ if $k\ge 0$, 
and zero otherwise. The Frobenius formula \cite{23} asserts
that
$$
\chi^\lam =
\left|
\matrix
\format\,\,\l&\quad\l&\quad\l&\quad\l&\quad\c&\quad\c\,\,\\
\eta^{\lam_1} & \eta^{\lam_1+1} & \eta^{\lam_1+2} & \hdots & \hdots & \hdots \\
\eta^{\lam_2-1} & \eta^{\lam_2} & \eta^{\lam_2 + 1} & \eta^{\lam_2+2} & 
\hdots & \hdots \\
\vdots & \eta^{\lam_3-1} & \eta^{\lam_3} & \eta^{\lam_3+1} & \hdots & \hdots \\
\vdots &  &  &\vdots  &  & \vdots \\
\vdots & \hdots &\hdots  & \eta^{\lam_i-i+j} & \hdots & \vdots \\
\vdots &  &  &\vdots  &  & \vdots
\endmatrix
\right|^\circ
$$
where the multiplication is to be understood as the external
one. The Frobenius formula expresses the character $\chi^\lam$ of the
irreducible representation of the symmetric group as a linear
combination of the functions $\eta^\lam$
$$
\eta^\lam=\eta^{\lam_1}\circ\eta^{\lam_2}\circ \dots \,.
$$
The following lemma is straightforward. 

\proclaim{Lemma 3}
A function $\phi\in \gZ(S(n))$ is multiplicative in the cycles of
a permutation if and only if for all
$n_1$, $n_2$ such that $n_1+n_2=n$, and for all functions
$f_1\in \gZ(S(n_1))$, $f_2\in \gZ(S(n_2))$ we have 
$$
\gathered
\la \phi|_{S(n_1)\times S(n_2)},f_1 \otimes f_2 \ra
_{S(n_1)\times S(n_2)} = \\ =
\la \phi|_{S(n_1)},f_1 \ra_{S(n_1)}
\la \phi|_{S(n_2)},f_2 \ra_{S(n_2)} \,.
\endgathered
\tag 2.6
$$
\endproclaim

\noindent{\bf Remark.}
The equality (2.6) is antilinear in $f_1$ and $f_2$. Therefore,  
it sufficient to check it for $f_1$ and $f_2$
ranging
independently over some linear bases of the spaces $\gZ(S(n_1))$ and
$\gZ(S(n_2))$.

\proclaim{Proposition 4}
The coefficients $m(\lam)$ correspond to a multiplicative central
function on the group $\Si$ if and only if we have 
$$
m(\lam)=\det\big[m(\lam_i-i+j)\big]_{i,j}
\tag 2.7
$$
for any Young diagram $\lam$. Here $m(k)$, $k\ge 0$, stands for
$m((k))$ and we set $m(k)=0$ for $k<0$.  
\endproclaim

\demo{Proof}
Let $\phi\in \gZ(\Si)$ be a multiplicative function. By the Frobenius
reciprocity,
$$
\la \phi|_{S(|\lam|)}, \eta^\lam \ra_{S(|\lam|)}=
\la \phi|_{S(\lam_1)\times \dots \times S(\lam_s)},1 \ra
_{S(\lam_1)\times \dots \times S(\lam_s)}\,.
$$
By the multiplicativity of $\phi$ this equals
$$
\prod_i \la \phi|_{S(\lam_i)},1 \ra _{S(\lam_i)}
=\prod_i m(\lam_i) \,.
$$
Therefore, by the Frobenius formula,
$$
m(\lam)=\la \phi|_{S(|\lam|)},\chi^\lam \ra_{S(|\lam|)}=
\det \big[m(\lam_i-i+j)\big]   \,.
$$

In the opposite direction, assume (2.7). Fix some $n$ and consider
the following dominance ordering of all partitions of $n$ 
$$
\mu \le \lam \quad
\Leftrightarrow \quad
\forall i \quad \mu_1+\dots+\mu_i \le
\lam_1+\dots+\lam_i \,. 
$$
Then (it is clear from
Frobenius formula) the transition matrix from the basis
$\{\chi^\lam\}$ to the basis $\{\eta^\lam\}$ is triangular.
The equality
$$
\la \phi, \eta^\lam \ra = \prod_i
{m(\lam_i)}
\tag 2.8
$$
is evident if $\lam=(n)$. Assume that
the equality is proved for all partitions $\mu$ of $n$, such that
$\mu > \lam$. Then it follows from (2.7) and Frobenius formula
that it is also true for $\lam$.
Therefore, the equality (2.8) is true  for all partitions
$\lam$.

Now let $\lam,\mu$ be two arbitrary partitions. Set $|\lam|=n_1$,
$|\mu|=n_2$. Denote by $\lam\cup\mu$ the union of parts of $\lam$
and $\mu$. Then
$$
\multline
\la \phi|_{S(n_1)\times S(n_2)},
\eta^\lam \otimes \eta^\mu \ra
_{S(n_1)\times S(n_2)} = 
\la \phi|_{S(n_1+n_2)},\eta^{\lam\cup\mu} \ra_{S(n_1+n_2)}\\
=\prod_i
{m(\lam_i)}
\prod_j
{m(\mu_j)}=\la \phi|_{S(n_1)},\eta^\lam \ra_{S(n_1)}
\la \phi|_{S(n_2)},\eta^\mu \ra_{S(n_2)}\,. 
\endmultline
$$
Since the functions $\eta^\lam$ form a basis in the
space of central functions, the 
multiplicativity of  $\phi$ is established.
\qed\enddemo

Recall that a sequence of real numbers $\{a_i\}$, $i=0,1,2,\dots$
is said to be {\it totally positive} if all the minors of the
following infinite Toeplitz matrix
$$
\left[
\matrix
a_0 & a_1 & a_2 & a_3 & \hdots \\
0 & a_0 & a_1 & a_2 & \ddots \\
 & 0 & a_0 & a_1& \ddots \\
 & & 0 & a_0 & \ddots \\
 & & & \ddots & \ddots
\endmatrix
\right]
$$
are non-negative. In particular, all the numbers $a_i$ are
non-negative.

A shift of indices and the multiplication of the entries by a positive
number preserve the total positivity, hence we
can always assume that $a_0=1$. Assume that for some $i>0$ we
have $a_i=0$. Then it follows from the inequality
$$
\det
\left[
\matrix
a_{i} & a_{i+1} \\
a_{0} & a_{1}
\endmatrix
\right] \ge 0
$$
that $a_j = 0$ for all $j \ge i$. Hence either all entries of a
totally positive sequence are positive or there exist such $N$ that
$$
a_i=0 \Leftrightarrow i\ge N\,.
$$
One easily checks that the 
multiplication by a positive geometric series preserves total
positivity. Therefore, if a sequence $\{a_i\}$ has at least two
positive entries, we can assume that
$$
a_0=a_1=1 \,.
$$
We call such a sequence {\it normalized totally positive}.

\proclaim{Proposition 5}
The map
$$
\phi \to \{m(i)\},\quad i=0,1,\dots
$$
is a bijection of the set of the characters of the group $\Si$ onto the
set of normalized totally positive sequences.
\endproclaim

\demo{Proof}
Let $\phi$ be a character of the group $\Si$. Clearly,
$m(0)=m(1)=1$. The inner product  of  $\phi|_{S(n)}$ 
with the trace of any
representation of $S(n)$ is always non-negative. In particular,
this is true for the trace $\chi^{\lam/\mu}$ of the
representation corresponding to the skew Young diagram
$\lam/\mu$. By multiplicativity of $\phi$ and
the analog of the Frobenius formula for $\chi^{\lam/\mu}$
\cite{28, \S1.5, (5.4)},
$$
\la \phi|_{S(|\lam/\mu|)},\chi^{\lam/\mu}
\ra_{S(|\lam/\mu|)} =
\det \big[ m(\lam_i-\mu_j-i+j) \big]_{ij} \ge 0 \,.
$$
One can easily see that all the minors involved in the definition
of total positivity can be obtained in this way.

In the opposite direction, assume that we are given a normalized
totally positive
sequence $\{m(i)\}$, $i=0,1,\dots$. Define a central function
$\phi_n$ on every symmetric group $S(n)$ by the formula 
$$
\la \phi_{n}, \eta^\lam \ra = \prod_i
{m(\lam_i)} \,, \quad |\lam|=n\,. 
$$
Because of $m(1)=1$ we have 
$$
\la \phi_{n+1}|_{S(n)}, \eta^\lam \ra =
\la \phi_{n+1}, \eta^{\lam\cup 1} \ra =
m(1) \prod_i
{m(\lam_i)} =  \la \phi_{n}, \eta^\lam \ra \,.
$$
This yields the coherence of $\{\phi_n\}$
$$
\phi_{n+1}|_{S(n)}=\phi_n \,.
$$
The multiplicativity of $\phi$ is clear, its positive definiteness
follows from the Frobenius formula, and $m(0)=1$ implies $\phi(1)=1$.
\qed\enddemo

It is convenient to form a generating series
$$
H(t)=1+m(1)\,t+m(2)\,t^2+\dots
$$
for  the numbers $m(i)$, $i=0,1,\dots$. As we
already know, either $H(t)$ is a polynomial, or all the terms of
$H(t)$ are non-zero. The series $H(t)$ has the following expression 
in terms of the numbers $c_i$, $i=1,2,\dots$. Recall that 
$$
c_k=\phi((12\dots k)) \,.
$$
Suppose that $\lam=1^{\rho_1} 2^{\rho_2} \dots$. We have 
$$
\align
m(k)&=\la \phi|_{S(k)},1 \ra_{S(k)} \\
&=\frac1{k!}\sum_{\lam,|\lam|=k}
{
\frac{k!}{\prod_i{\rho_i!i^{\rho_i}}}
\prod_i{c_i^{\rho_i}}
} \\
&=\sum_{\lam,|\lam|=k}
{\prod_i
{\frac 1{\rho_i!}
\left(\frac{c_i}{i}\right)^{\rho_i} }}  \quad .
\endalign
$$
Therefore,
$$
\align
H(t)&=\sum_{k=0}^\infty {t^k \sum_{\rho_1+2\rho_2+3\rho_3\dots=k}
{\prod_i
{\frac 1{\rho_i!}
\left(\frac{c_i}{i}\right)^{\rho_i} }}} \\
&=\prod_{i=1}^\infty
\sum_{j=0}^\infty
{\frac 1{j!}
\left(\frac{c_i t^i}{i}\right)^{j} } \\
&=\prod_{i=1}^\infty
{\exp
\left(\frac{c_i t^i}{i}\right) } =
\exp
\left(\sum_{i=1}^\infty \frac{c_i}{i} t^i\right)   \quad .
\endalign
$$
Vice versa, the numbers $c_i$ can be easily obtained from the
series $H(t)$,
$$
c_1+c_2 t +c_3 t^2 + \dots = \frac{H'(t)}{H(t)} \,.
$$
We shall state the Thoma theorem in the following form.

\proclaim{Theorem \rm (Thoma \cite{49})}
The generating functions for Fourier coefficients $\{m(i)\}$,
$i=0,1,\dots$ of the characters of the group $\Si$ have the form
$$
H(t)=e^{\gamma t} \prod_i \frac{1+\be_i t}{1-\a_i t}
$$
where
$$
\a_i \ge 0,\quad \be_i \ge 0,\quad \gamma \ge 0,
\quad \sum \a_i+ \sum \be_i + \gamma = 1\,.
$$
\endproclaim

This is equivalent to the following: 

\proclaim{Theorem \rm (Edrei \cite{17})}
The generating functions of totally positive sequences
have the form
$$
H(t)=e^{\gamma t} \prod_i \frac{1+\be_i t}{1-\a_i t}
$$
where
$$
\a_i \ge 0,\quad \be_i \ge 0,\quad \gamma \ge 0,
\quad \sum \a_i+ \sum \be_i   < \infty \,.
$$
\endproclaim

\demo{Proof of the Thoma theorem}
First, assume that $H(t)$ is not a polynomial. Then it follows from the
inequalities
$$
\det
\left[
\matrix
m(i) & m(i+1) \\
m(i-1) & m(i)
\endmatrix
\right] \ge 0
$$
that
$$
1=
\frac{m(1)}{m(0)}\ge
\frac{m(2)}{m(1)}\ge
\frac{m(3)}{m(2)}\ge
\dots \ge 0 \,.
\tag 2.9
$$
Therefore, there exists a limit
$$
\a=\lim \frac{m(n+1)}{m(n)}, \quad n\to\infty \,.
$$
Obviously, $0\le\a\le 1$ and $1/\a$ is the
convergence radius of the series $H(t)$.

If $\a=1$  then by (2.9) we have $m(i)\equiv 1$ and
$$
H(t)=\frac 1{1-t} \,.
$$
This $H$-series corresponds to the trivial representation. We
denote it by $H_{1}(t)$.

Assume now that $0<\a<1$. Then we construct, starting from the
character $\phi$, a new character $\pht$ with Fourier
coefficients
$$
\mt(\lam)=(1-\a)^{-|\lam|}\lim \frac{m(n\cup\lam)}{m(n)},
\quad n \to \infty \,.
$$
It is easy to compute these limits explicitly. One has
$$
\align
&\lim \frac{m(n\cup\lam)}{m(n)}=
\lim
\frac {1}{m(n)}
\det
\left[
\matrix
m(n) & m(n+1) & m(n+2) & \hdots \\
m(\lam_1-1) & m(\lam_1) & m(\lam_1+1) & \hdots \\
m(\lam_2-2) & m(\lam_2-1) & m(\lam_2) &  \hdots \\
\vdots &  \vdots  &  \vdots & \ddots
\endmatrix
\right] \\
&=
\lim
\det
\left[
\matrix
1 & m(n+1)m(n)^{-1} & m(n+2)m(n)^{-1} & \hdots \\
m(\lam_1-1) & m(\lam_1) & m(\lam_1+1) & \hdots \\
m(\lam_2-2) & m(\lam_2-1) & m(\lam_2) &  \hdots \\
\vdots &  \vdots  &  \vdots & \ddots
\endmatrix
\right] \\
&=
\det
\left[
\matrix
1 & \a & \a^2 & \hdots \\
m(\lam_1-1) & m(\lam_1) & m(\lam_1+1) & \hdots \\
m(\lam_2-2) & m(\lam_2-1) & m(\lam_2) &  \hdots \\
\vdots &  \vdots  &  \vdots & \ddots
\endmatrix
\right]  \\
&=
\det
\left[
\matrix
1 & 0 & 0 & \hdots \\
m(\lam_1-1) & m(\lam_1)-\a m(\lam_1-1) & m(\lam_1+1)-\a m(\lam_1) & \hdots \\
m(\lam_2-2) & m(\lam_2-1)-\a m(\lam_2-2) & m(\lam_2)-\a m(\lam_2-1) &  \hdots \\
\vdots &  \vdots  &  \vdots & \ddots
\endmatrix
\right]  \\
&=
\det
\left[
\matrix
m(\lam_1)-\a m(\lam_1-1) & m(\lam_1+1)-\a m(\lam_1) & \hdots \\
m(\lam_2-1)-\a m(\lam_2-2) & m(\lam_2)-\a m(\lam_2-1) &  \hdots \\
\vdots  &  \vdots & \ddots
\endmatrix
\right]\quad .
\endalign
$$
In particular,
$$
\mt(k)=(1-\a)^{-k}(m(k)-\a m(k-1)) \,.
\tag 2.9
$$
Therefore,
$$
\mt(\lam)=\det\big[\mt(\lam_i-i+j)\big]_{i,j}\,.
$$
We see that the Fourier coefficients $\{\mt(\lam)\}$ determine a
multiplicative positive definite function on the group $\Si$. Set
$$
\Ht(t)=\sum \mt(k)\, t^k \,.
$$
Then the equality (2.9) means that
$$
H(t)=H_{1}(\a t)\,\Ht((1-\a)t)\,.
$$
The representation-theoretical meaning of this factorization will
be made clear in Section 3. Note that the convergence radius of
the series $\Ht((1-\a)t)$ is not less than that of the series
$H(t)$. Iterating this factorization, we obtain the following factorization
$$
H(t)=
H_{\text{ent}}\left(\left(1-\sum \a_i\right) t\right) \, 
\prod_i H_{1} (\a_i t)
 \,,
$$
where
$$
\a_i \ge 0, \quad \sum \a_i \le 1 \,,
$$
and $H_{\text{ent}}$ is an entire function which is also an $H$-series
for a character of the group $\Si$.

Now, instead of $\phi$,  consider the character $\phi^-$,
$$
\phi^-(g)=\sgn(g)\phi(g),\quad g\in \Si \,.
$$
In particular,
$$
c^-_k=(-1)^{k-1}c_k \,.
$$
Therefore, the character $\phi^-$ is associated with the
$H$-series
$$
H^-(t)=
\exp \left(\sum_{i=1}^\infty (-1)^{i-1} \frac{c_i}{i} t^i\right) =
\exp \left(- \sum_{i=1}^\infty \frac{c_i}{i} (-t)^i\right) =
H(-t)^{-1}\,.
$$
Now factoring the function $H_{\text{ent}}(-t)^{-1}$ we obtain
$$
H(t)=
H^\ast(\gamma t)\,
\prod H_{1} (\a_i t) \,
\prod H_{\text{sgn}} (\be_i t) \,
 \quad \sum\a_i+\sum\be_i+\gamma=1 \,,
$$
where
$$
H_{\text{sgn}}(t)=1+t \,,
$$
and $H^\ast(t)$ is an $H$-series for a character of the group
$\Si$ which is an entire function with no zeros.
Hence the series
$$
c_1+c_2 t +c_3 t^2 + \dots = \frac{H'(t)}{H(t)}
$$
represents an entire function.  Now, and this is
the whole point of this section, we shall show that this is only
possible if
$$
c_1=1,\quad c_2=c_3= \dots = 0\,.
$$
Indeed, we have
$$
c_k=\int t^{k-1}d\mu\,.
$$
Assume that the measure $\mu$ is not concentrated at $t=0$.
Then for some $\e>0$ we have
$\mu([-1,1]\setminus[-\e,\e])=const>0$, hence
$$
c_{2k+1}=\int t^{2k}d\mu\ge const\,\,\e^{2k}
$$
for all $k$, and
hence the convergence radius of the series $\sum c_kz^{k-1}$ does
not exceed $1/\e$. Therefore, the measure $\mu$ should be
supported at zero, which implies  $c_2=c_3=\dots=0$. In other words,
$$
H^\ast(t)=e^t \,. \qed
$$
\enddemo

\noindent{\bf Remark.}
The case of entire function without zeros and poles was actually
the most difficult part in the original proof \cite{17,49}.

\subhead 2.5 Classification of the irreducible admissible
representations of the pair $\GKd$ \endsubhead

Let $\pi$ be an irreducible admissible representation of depth
$d$ determined by a Thoma measure $\mu$ and  Young
distributions $\L$, $\Mu$.

Denote by $s$ and $t$ the permutations
$(1,d+1)$, $(-1,-d-1)\in\Gd$. 

\proclaim{Lemma 4} Suppose $x\in\supp\L$ and suppose a vector
$\ze\in R(\pi)$, $\|\zeta\|=1$, satisfies $A_1\ze=x\ze$. 
Set $\ze^{(x)}=\dt_x(A_{d+1})\ops \ze$. Then 
$$
\alignat2
&\text{\rm{a)}}\quad&&(\ze^{(x)},\ze^{(x)})=\mu(x),\\
&\text{\rm{b)}}\quad&&(s\ze^{(x)},\ze^{(x)})=x,\\
&\text{\rm{c)}}\quad&&(st\ze^{(x)},\ze^{(x)})=0\,.\\
\intertext{If, in addition, $x\in\supp\Mu$ and $A_{-1}\ze=x\ze$ then}
&\text{\rm{d)}}\quad&&(t\ze^{(x)},\ze^{(x)})=x\,.\\
\endalignat
$$
\endproclaim

\demo{Proof}
Part a) follows from the definition of the spectral measure. By
Proposition 3, $\mu(x)\ne0$ if $x\ne 0$. Part b) has actually been
already established in the proof of Proposition 3. Part d) then
follows from the fact that, by virtue of (1.4),
$$
A_{-d-1}^k\ze=A_{d+1}^k\ze\,, \quad k=1,2,\dots \,,
$$
and hence
$$
\ze^{(x)}=\dt_x(A_{-d-1}) \ze \,.
$$
Let us check part c). The operator $P_d s t P_d$ is conjugate to
$P_{d-1}$, and therefore equals zero since $H(\pi)_{d-1}=0$ by
assumption. Hence,
$$
\split
&0=(P_d st P_d \ze,\ze)=(st\ze,\ze)=(st\dt_x(A_1)\ops \ze,\dt_x(A_1)\ops \ze)\\
&\overset(1.1),(1.2)\to=(st\dt_x(A_{d+1})\ops \ze,\dt_x(A_{d+1})\ops \ze)
=(st\ze^{(x)},\ze^{(x)})\,.\qed
\endsplit
$$
\enddemo

\proclaim{Theorem 3}
An irreducible admissible
representation of depth $d$ with the Thoma measure $\mu$ and
Young distributions $\L$, $\Mu$, $|\L|=|\Mu|=d$ exists if
and only if for all $x\in[-1,1]$ we have
$$
\alignat2
\ell(\L(x))+\ell(\Mu(x))&\le \mu(x)/|x|,&\qquad&x>0,
\\
\ell(\L'(x))+\ell(\Mu'(x))&\le \mu(x)/|x|,&\qquad&x<0,
\endalignat
$$
where the prime denotes the transposition   of Young
diagrams and $\ell$ stands for the number of rows in a Young 
diagram.
\endproclaim

\demo{Proof}
The sufficiency was proved in \cite{32} using an explicit
construction of the representation. We shall give the explicit
construction of  all corresponding representations in Section 3.

Let us prove the necessity. Assume that
$$
x\in\supp\L\cup\supp\Mu \quad  \text{ and } \quad  x>0 \,.
$$
Set
$$
l_1=\ell(\L(x)), \quad l_2=\ell(\Mu(x))\,.
$$
We shall assume that $l_1,l_2>0$.  The case where $l_1=0$ or
$l_2=0$ is more simple, and can be dealt with analogously. Denote by
$S(l_1)\times S(l_2)$ the subgroup in the group $\Gd(d)$
which permutes the numbers $\{1,\dots,l_1\}$ and $\{-1,\dots,-l_2\}$.

By the branching rule for the representations of finite
symmetric groups and the definition of $\L$ and $\Mu$, there
exists a vector $\ze\in R(\pi)$ which is anti-invariant with
respect to the group $S(l_1)\times S(l_2)$ and also satisfies
$$
A_i\ze=x\ze, \quad i=1,\dots,l_1,-1,\dots,-l_2\,.
$$
Then the vector $\ze^{(x)}$ is $S(l_1)\times
S(l_2)$-anti-invariant, too.

Let $\Alt$ be the operator of anti-symmetrization over the group
$S(l_1+1)\times S(l_2+1)$ permuting $\{1,2,\dots,l,d+1\}$ and
$\{-1,-2,..,-l,-d-1\}$. The function 
$$
(\sgn\ops (g)\ops g\ze^{(x)},\ze^{(x)})
$$
is invariant with respect to right and left translations by the
elements of the group $S(l_1)\times S(l_2)$. The group
$S(l_1+1)\times S(l_2+1)$ consists of four double cosets with
respect to subgroup $S(l_1)\times S(l_2)$; their representatives
are $1$, $s$, $t$, $st$, and their cardinalities are
$l_1!\ops l_2!$, $l_1l_1!\ops l_2!$, $l_2l_1!\ops l_2!$,
$l_1l_2l_1!\ops l_2!$. Hence, by  the above lemma we have
$$
\split
(\Alt\ze^{(x)},\ze^{(x)})
&=\frac{l_1!\ops l_2!\ops \mu(x)-l_1l_1!\ops l_2!\ops x-l_2l_1!\ops l_2!\ops x
+l_1l_2l_1!\ops l_2!\ops 0}{(l_1+1)!\ops (l_2+1)!}\\
&=\frac{x(\mu(x)/x-l_1-l_2)}{(l_1+1)(l_2+1)}\,.
\endsplit
$$
Since the operator $\Alt$ is a projection, the result is
non-negative, and therefore
$$
l_1+l_2\le \mu(x)/x\,.
$$
The case $x<0$ can be reduced to that of $x>0$ by replacing $T$
with $T\otimes(\sgn\otimes\sgn)$.
\qed\enddemo

\noindent{\bf Remark.}
From the point of view of the figure in Section 1.2, the use of spectral
projectors $\dt_x(A_i)$ allows one to get rid of  all Young
diagrams but those growing from the point $x$. In Section 3, we
shall be concerned with the opposite problem: how to plant a
Young diagram at a given point of the interval $[-1,1]$.

\subhead  2.6 Description of  $\Ke\backslash\Ge/\Ke$ \endsubhead

In this section  $G$ stands the group $\Ge(n)$
and  $K$  denotes the subgroup $\Ke(n)$. We shall
recall basic facts about the cosets $\KGK$ which we shall need
in the next section. Denote by $S(n)\subset
\Ge(n)$ the subgroup fixing the points $\{-1,\dots,-n\}$.

The set $\KGK$ clearly coincides with the set of orbits of the
group $K$ on the set $G/K$. The set $G/K$ is naturally identified
with the set $\Pi$ of partitions of the set $\{\pm 1,\pm
2,\dots,\pm n\}$ into pairs (note the difference
between partitions of a set and partitions of a number).
The group $K$ is the stabilizer
of the partition $\sigma$
$$
\sigma=\left\{
\left\{\pm 1\right\},\dots,
\left\{\pm n\right\} \right\}\,.
$$
For any two partitions $\tau$, $\upsilon$ let $\tau \lor
\upsilon$ denote their least upper bound, that is, the
finest partition consisting of whole blocks of $\tau$
and $\upsilon$. If $\tau$ and $\upsilon$ were
partitions into  pairs (more generally, into even  blocks), 
then the partition $\tau \lor \upsilon$ is also a
partition into  even blocks. Therefore the block cardinalities of
$\tau\lor\upsilon$, divided by two, form a partition of $n$ which
we denote by $\tau\tri\upsilon$. The function
$f(\tau)=\tau\tri\sigma$ with values in the set of
partitions of $n$ is an invariant of the action of the group $K$
on the space $\Pi$. The following proposition is well known
\cite{46, 32}.

\proclaim{Proposition}
The function $f(\tau)=\tau \tri \sigma$ separates the orbits of
the group $K$ in $\Pi$.
\endproclaim

\proclaim{Corollary}
\roster
\item The set $\KGK$ is parameterized by partitions $\lam$ of the
number $n$.
\item The intersection of the double coset corresponding to a
partition $\lam$ with the subgroup $S(n)$ consists of
permutations with the cycle structure $\lam$. In particular,
this intersection is non-empty.
\endroster
\endproclaim

Let $\lam=1^{m_1}2^{m_2}\dots$ be the partition with $m_i$ parts
of size $i$. Denote by $\ell(\lam)$ the number of parts in the
partition $\lam$. Set $z_\lam = \prod{i^{m_i} m_i!\,}$. Denote by
$K\lam K$ the double coset corresponding to the partition $\lam$.
The following proposition can be established by a direct
combinatorial argument \cite{46}:

\proclaim{Proposition 7}
$$
|K\lam K|=2^{2n-\ell(\lam)}\frac{(n!)^2}{z_\lam}
$$
\endproclaim

\proclaim{Corollary}
For an element $g \in \Ge(n)$, let $\ell(g)$ be the number of
parts in the partition corresponding to the double coset $KgK$.
If $t$ is a formal variable, then
$$
\sum_{g \in \Ge(n)} {t^{\ell(g)}}=
n! \, 2^n t (t+2) (t+4) \cdots (t + 2n-2)\,.
$$
\endproclaim

\demo{Proof}
We have
$$
\split
\sum_{g \in \Ge(n)} {t^{\ell(g)}}
&=\sum_{\lam \vdash n}
2^{2n-\ell(\lam)}(n!)^2 z_\lam^{-1} t^{\ell(\lam)}\\
&=n! \, 2^{2n} \sum_{\lam \vdash n}
n! \, z_\lam^{-1} \left(\frac{t}{2}\right)^{\ell(\lam)}
\intertext{
By the identity (2.4) for the Stirling numbers
\cite{45} already employed in the proof of Theorem 2, this is}
&=n! \, 2^{2n} \frac t2 \left(\frac t2+1\right)
\left(\frac t2 +2\right) \cdots \left(\frac t2 + n-1\right)\\
&=n! \, 2^n t (t+2) (t+4) \cdots (t + 2n-2)\,.\qed
\endsplit
$$
\enddemo

\subhead 2.7 Spherical representations of the pair $\GKe$ \endsubhead

Let $\pi$ be a spherical representation of the pair $\GKe$ and let
$\phi$ be the corresponding spherical function. We know from the
previous section that the
group $\Ge$ is the product of its subgroups $\Ke \Si \Ke$. Hence
$\phi$, as a $\Ke$-biinvariant function, is completely determined
by its restriction to the subgroup $\Si$. This restriction is a
normalized positive definite function. Just as in the case of
spherical representations of the pair $\GKd$,  one checks that
$$
\phi(\si)=\prod_{k\in [\si]} c_k \,, \quad \si\in\Si\,.
$$
That is,  the function $\phi$ is multiplicative in the cycles of a
permutation and,  hence, by Proposition 2,  $\phi$ is
indecomposable. Therefore, it has
the form $\phm$, for some Thoma measure $\mu$. Since
for any $g \in \Si$ the intersection
$$
\Si \cap \Ke g \Ke
$$
is the  conjugacy class of $g$ in the group $\Si$, every function $\phm$
has a unique $\Ke$-biinvariant
extension to the group $\Ge$. We denote this extension by
$\tphm$.

To summarize, the spherical functions of the pair $\GKe$ are
precisely those functions $\tphm$ that are positive
definite on $\Ge$. The description of this set is given by the
following

\proclaim{Theorem 4}
The function $\tphm$ is a spherical function of the pair $\GKe$
if and only if
$$
\frac{\mu(x)}{|x|}\in 2\bZp
$$
for every $x\in [-1,0)$.
\endproclaim

\demo{Proof}
The sufficiency of the condition follows from the explicit
constructions of representations (see Chapter 3). Let us prove
the necessity.

Let $\pi$ be the spherical representation of the pair $\GKe$
corresponding to the spherical function $\tphm$. Let $\xi$ be the
spherical vector of the representation $\pi$. Take some $x \in
[-1,0)$. In the same way as we already did it in the proof
of Theorem 2, by replacing the vector $\xi$ by the vectors
$$
\xi^{(m)}=
\prod_{i=1}^m {\delta_x(A_i)} \xi\,, \quad m=1,2,3,\dots\,,
$$
one can effectively reduce the case of a general Thoma
measure $\mu$ to that of Thoma measure supported at a single point $x$.
Thus, we can assume that
$$
\supp \mu = \{x\}, \quad \mu(x)=1\,.
$$
In this case the value of the function $\tphm$ at the element
$g\in \Ge(n)$ is
$$
\tphm(g)=x^{n-\ell(g)}\,.
$$
Consider the non-negative expression
$$
\sum_{g \in \Ge(n)} {\tphm(g)}
= \sum_{g \in \Ge(n)} {x^{n-\ell(g)}}\,.
$$
By the Corollary to Proposition 7 it equals
$$
\gathered
n! \, 2^{n} x^n \frac 1x \left(\frac 1x+2\right)
\left(\frac 1x +4\right) \cdots \left(\frac 1x + 2n-2\right) =\\=
n! \, 2^n 1 (1+2x) (1+4x) \cdots (1 + (2n-2)x)\,.
\endgathered
$$
Since $x<0$, this product can be non-negative for all $n$ only if
it terminates, which happens  if
$$
x=-\frac 1{2k}
$$
for some $k\in \bN$.  That is,
$$
\nu(x)= 2k\,. \qed
$$
\enddemo

This result is parallel to the following well known result
from representation theory of finite symmetric groups. Namely, an
irreducible admissible representation of the group $\Ge(n)$
corresponding to a Young diagram $\lam$, $|\lam|=2n$, contains a
$\Ke(n)$-invariant vector if and only if all the parts of $\lam$ are even
\cite{28, \S1.8, Example 6}.

\subhead  2.8 Classification of irreducible admissible
representations of the pairs $\GKo$, $\GKe$ \endsubhead

Let $\pi$ be an irreducible admissible representation of depth
$d$ of  the Gelfand pair $\GKe$. The case of the pair $\GKo$ is entirely
analogous. The representation $\pi$ is determined by a Young
distribution $\L$, such that $|\L|=2d$, and  a Thoma
measure $\mu$. By the results of the previous section,
$\mu(x)/|x|\in2\bZ_+$ for all $x<0$.

The proof of the next proposition is a word for word copy of that
of Proposition 3 in Section 2.2.

\proclaim{Proposition 8}
$\supp\L\subset\supp \mu\cup\{0\}$\,.
\endproclaim

The description of irreducible admissible representations of the
pair $\GKe$ is provided by the following theorem.

\proclaim{Theorem 5}
An irreducible admissible representation if the pair $\GKe$ (respectively,  of
the pair $\GKo$) with depth $d$, Thoma measure $\mu$, and  Young
distribution $\L$, where $|\L|=2d$ (resp.\ $2d+1$), exists if and only if
$\mu(x)/|x|\in2\bZ_+$ for all $x<0$ and
$$
\alignat3
\L'(x)_1+\L'(x)_2&\le\, &&\mu(x)/|x|,&\qquad&x>0,\\
\L(x)_1&\le\, &&\mu(x)/2|x|,&\qquad&x<0
\endalignat
$$
for all $x\in[-1,1]$.
\endproclaim

\demo{Proof}
The sufficiency follows from the explicit construction of
representations \cite{32} to be discussed in Section 3.

Let us prove the necessity. Suppose $x>0$ and set
$$
l_1=\L'(x)_1, \quad l_2=\L'(x)_2 \,.
$$
Consider the main case when $l_1,l_2>0$. Other cases are
similar.  Denote by $S(l_1)\times
S(l_2)$ the subgroup in the group $\Ge(d)$ which permutes the numbers
$\{d-l_1+1,\dots,d\}$ and $\{-d,\dots,-d+l_2-1\}$.

By the branching rule for representations of finite symmetric
groups and the definition of $\L$, there exists a vector
$\ze\in R(\pi)$ which is antiinvariant under the action of the
group $S(l_1)\times S(l_2)$ and satisfies 
$$
A_i\ze=x\ze, \quad i=d-l_1+1,\dots,d,-d,\dots,-d+l_2-1\,.
$$
Then the vector
$$
\ze^{(x)}=\dt_x(A_{d+1})\ze
$$
is also $S(l_1)\times S(l_2)$-antiinvariant.

Let $\Alt$ be the antisymmetrization over the group
$S(l_1+1)\times S(l_2+1)$ which permutes the numbers
$\{d-l_1+1,\dots,d,d+1\}$ and $\{-d-1,-d,\dots,-d+l_2-1\}$. The
function 
$$
(\sgn\ops (g)\ops g\ze^{(x)},\ze^{(x)})
$$
is $S(l_1)\times S(l_2)$-biinvariant. Denote by $s$ and $t$ the
permutations $(d,d+1)$, $(-d,-d-1)\in \Ge$. Arguing as in the
proof of Lemma 4 we obtain
$$
\aligned
&(\ze^{(x)},\ze^{(x)})=\mu(x)\\
&(s\ze^{(x)},\ze^{(x)})=x \\
&(t\ze^{(x)},\ze^{(x)})=x\\
&(st\ze^{(x)},\ze^{(x)})=0 \,.
\endaligned \tag 2.10
$$
The group $S(l_1+1)\times S(l_2+1)$ consists of four double
cosets of the subgroup $S(l_1)\times S(l_2)$; their
representatives are $1$, $s$, $t$, $st$ and the corresponding
cardinalities are $l_1!\ops l_2!$, $l_1l_1!\ops l_2!$, $l_2l_1!\ops l_2!$,
$l_1l_2l_1!\ops l_2!$. Therefore, by virtue of (2.10),
$$
\split
(\Alt\ze^{(x)},\ze^{(x)})
&=\frac{l_1!\ops l_2!\ops \mu(x)-l_1l_1!\ops l_2!\ops x-l_2l_1!\ops l_2!\ops x
+l_1l_2l_1!\ops l_2!\ops 0}{(l_1+1)!\ops (l_2+1)!}\\
&=\frac{x(\mu(x)/x-l_1-l_2)}{(l_1+1)(l_2+1)}\,.
\endsplit
$$
Since the result should be
non-negative, we conclude that
$$
l_1+l_2\le \mu(x)/x\,.
$$
Now suppose that  $x<0$ and set
$$
l=\L(x)_1 \,.
$$
Denote by $S(l)$ the subgroup in the group $\Ge(d)$ which permutes the
numbers $\{d-l+1,\dots,d\}$.

Again, there exists a vector $\ze\in R(\pi)$ invariant under the action
of the group $S(l)$ and such that
$$
A_i\ze=x\ze, \quad i=d-l+1,\dots,d \,.
$$
Then the vector
$$
\ze^{(x)}=\dt_x(A_{d+1})\ze
$$
is also $S(l)$-invariant. The permutation $r=(d+1,-d-1)$ belongs
to the group $\Ke_d$. By the definition of the subspace $R(\pi)$,  the
vector $\ze\in R(\pi)$ will be invariant under $r$,
hence with respect to the group $S(l)\times
S(2)$ which permutes the numbers $\{d-l+1,\dots,d\}$ and
$\{d+1,-d-1\}$. Let $\Sym$ be the symmetrization over
the group $S(l+2)$ which permutes the numbers
$\{d-l+1,\dots,d,d+1,-d-1\}$. Denote by $s$ and $t$ the
permutations $(d,d+1)$, $(d-l+1,-d-1)\in \Ge$. As usual,
$$
\align
&(\ze^{(x)},\ze^{(x)})=\mu(x)\\
&(s\ze^{(x)},\ze^{(x)})=x \\
&(st\ze^{(x)},\ze^{(x)})=0 \,.
\endalign
$$
The group $S(l+2)$ consists of three double cosets of
 $S(l)\times S(l_2)$; their representatives are $1$, $s$,
$st$ and the cardinalities are $2l!$, $4l\,l!$, $l(l-1)l!$.
Therefore,
$$
\split
(\Sym\ze^{(x)},\ze^{(x)})
&=\frac{2l!\ops \mu(x)+4l\,l!\ops x
+l(l-1)l!\ops 0}{(l+2)!}\\
&=\frac{2|x|(\mu(x)/|x|-2l)}{(l+1)(l+2)}\,.
\endsplit
$$
The result should be
non-negative, whence
$$
l \le \,\mu(x)/2|x|\,. \qed
$$
\enddemo

\noindent{\bf Remark.}
The cases $x>0$ and $x<0$ for the pairs $\GKe$ and $\GKo$ are not
symmetric. This is because the function $\sgn$ on
the group $\Si$ cannot be extended to a $\Ke$-invariant positive
definite function on the group $\Ge$.

\head 3.~Construction of representations \endhead

The main object considered in this section 
is a certain operation on admissible representations which
we call {\it mixing} the representations. This construction
is very much parallel to the Olshanski's construction in \cite{32}
and only slightly more general (see also \cite{6,56}). It yields an explicit
construction of actually all irreducible admissible representations
whereas the methods of \cite{32} produce only an open
subset in of the admissible dual. The mixture of
representations is, essentially, a special sort of
an induced representation  as we shall see in Section 3.2.

In Sections 3.1--3.3 we shall deal with the pairs $\Ge$ and $\Go$
(mainly with $\Ge$). We shall comment briefly on the case of
$\Gd$ in Section 3.4. 

\subhead 3.1  Mixtures of representations \endsubhead

Let $\pi_1$ and $\pi_2$ be two admissible representations of the
pair $\GKe$. Let $p_1$ and $p_2$ be two numbers, such that
$p_1>0$, $p_2>0$, $p_1+p_2=1$. We shall define the {\it mixture} of
representations $\pi_1$ and $\pi_2$ with the weights $p_1$ and
$p_2$. One can similarly define the mixture of admissible representation of
the pair $\GKe$ with an admissible representation of the pair
$\GKo$ or the mixture of two admissible representations of the
pair $\GKo$. Set $d_i=d(\pi_i)$, $i=1,2$.

In the set of all functions $f:\bZz \to \{1,2\}$ consider the
following subset $X$,
$$
X=\left\{ f|
f\left( i \right)=f\left( -i \right)
\text{ for almost all } i
\right\} \,,
$$
where ``almost all'' means ``all but finitely many''. The set
$X$ is a union of an increasing sequence of subsets
$$
X_0 \subset X_1 \subset  X_2 \subset \dots
$$
where
$$
X_n=\left\{ f\in X|
f\left( i \right)=f\left( -i \right)
,|i| > n
\right\}   \,.
$$
The map
$$
X_n\owns f\mapsto (f(-n),\dots,f(-1),f(1),\dots)\in \{1,2\}^\infty
$$
is a bijection. We transfer the product topology from
$\{1,2\}^\infty$
to $X_n$ via this map and endow the set  $X$ with the direct limit topology.
Consider the following  measure $\ompp$ on $X$,
$$
\ompp\left( \left\{
f|f\left( i \right)=f_i,|i|\leq n,
f\left( i \right)=f\left( -i \right),
i>n
 \right\} \right)=
\prod_{|i|<n}{p_{f_i}^{1/2}} \,.
$$
On each set $X_n$ the measure $\ompp$ is finite. The group $\Ge$
acts on $X$ and preserves the measure $\ompp$.

It follows from the definition of the space $X$ that for every
 $f\in X$ the parity of the number $|f^{-1}(1)\cap\{-N,\dots,N\}|$ 
stabilizes as $N\to\infty$. Consider the subset $Y \subset X$,
$$
Y=\left\{f,
|f^{-1}(1)\cap \{-N,\dots,N\}|\in 2\bZ, N\gg 0, 
|f^{-1}(i)|=\infty,i=1,2
 \right\}   \,.
$$
This subset is measurable and $\ompp(Y)>0$, since the condition
$|f^{-1}(i)|=\infty$, $i=1,2$ means deletion of a countable
set of zero measure. It is also clear that this subset is
$\Ge$-invariant. We set $Y_n=Y\cap X_n$.

For every $f\in Y$ there exist unique bijections
$\eta_i:f^{-1}(i)\rightarrow\bZz$, such that
$\eta_i(-a)=-\eta_i(a)$ for almost all $a$ and $a<b$ if and only if
$\eta_i(a)<\eta_i(b)$. Define a cocycle
$$
c:\Ge \times Y \to \Ge \times \Ge
$$
on the generators $(i,i+1)$ of the group $\Ge$ by the formula
$$
c\left((i,i+1),f\right)=
\cases (e,e), &f(i)\neq f(i+1) \\
\left(\left(\eta_1(i),\eta_1(i+1)\right),e\right), &f(i)=f(i+1)=1 \\
\left(e,\left(\eta_2(i),\eta_2(i+1)\right)\right), &f(i)=f(i+1)=2  \,.
\endcases
$$
Denote by $H$ the Hilbert space of maps
$$
F:Y\rightarrow H\left( \pi_1 \right)
\otimes H\left( \pi_2 \right)
$$
with inner product
$$
(F_1,F_2)_H=\int_Y\left(F_1(f),F_2(f)\right)_
{H\left( \pi_1 \right)\otimes H\left( \pi_2 \right)}
\,\ompp(df)   \,.
$$
Define a representation of the group $\Ge$ in the space $H$ by
the formula
$$
\left[g\cdot F\right](f)=
\pi_1\otimes\pi_2
\left(c(g,g^{-1}\cdot f)\right) F (g^{-1}\cdot f)   \,.
$$
This representation is unitary.

Now, our next goal is to compute the subspaces of invariants
$H_n$. Fix two numbers $a,b \in \bZp$. Set
$K=2a+2b$ and define
$$
\split
\Yab=\left\{f \in Y \right.,
&|f^{-1}(1)\cap \{-K,\dots,K\}|=2a \\
&|f^{-1}(2)\cap \{-K,\dots,K\}|=2b \\
&\left. f(i)=f(-i),|i|>K \right\}    \,.
\endsplit
$$
Let $D$ range over all subsets of cardinality $2a$ of the set
$\{-K,\dots,K\}\setminus 0$. Set
$$
\YabD=\left\{f \in \Yab |
f^{-1}(1)\cap \{-K,\dots,K\}=D
\right\}  \,.
$$
Denote by $\supp F$ the complement to the largest open subset
where $F$ is equal to zero almost everywhere. Denote by $\HYab$
the space of maps $F$, such that $\supp F \subset \Yab$. Denote
by $\Hab \subset \HYab$ the subspace of maps with the support in
$\Yab$, which are constant on all of $\YabD$ and take the values
in $H\left( \pi_1 \right)_a\otimes H\left( \pi_2 \right)_b$.

\proclaim{Proposition 9}
$$
H_n=\bigoplus_{a+b=n}{\Hab}    \,.
$$
\endproclaim

\demo{Proof}
The inclusion $\bigoplus_{a+b=n}{\Hab}\subset H_n$ is obvious.
Let us check the inverse inclusion. Assume that $F\in H_n$, and
let the numbers $r$, $s$ be such that $s\geq r>n$. Set
$$
\Wrs=\left\{f \in Y |
f(r)\neq f(-r),f(i)=f(-i),|i|>s \right\}   \,.
$$
The sets $\Wrs$ are open. Let $l_1,l_2$ be two distinct integers, such
that $l_1,l_2>s$. Consider the action of permutations
$(r,l_i)(-r,-l_i)\in\Ke_n$ on the sets $\Wrs$. We have
$$
(r,l_1)(-r,-l_1)\cdot\Wrs \cap
(r,l_2)(-r,-l_2)\cdot\Wrs =\emptyset
$$
$$
\ompp \left((r,l_1)(-r,-l_1)\cdot\Wrs \right)=
\ompp \left((r,l_2)(-r,-l_2)\cdot\Wrs \right)      \,.
$$
Since the map $F$ is square summable, it should vanish almost
everywhere on each set $\Wrs$. This means that
$$
\supp F \subset Y \setminus
\bigcup_{n<r\leq s}{\Wrs} =
\bigcup_{a+b=n}{\Yab}    \,.
$$
In other words, $F \in \bigoplus_{a+b=n}{\HYab}$.

All the subspaces $\HYab$ are invariant with respect to $\Ke_n$.
The subspaces of invariants are always consistent with
decompositions in a direct sum, hence
$$
F \in \bigoplus_{a+b=n}{\left(\HYab\right)_n}   \,.
$$
Let $F\in\left(\HYab\right)_n$ for some $a,b$. For every
$g\in\Ke_n$ the equality $\left[g\cdot F\right](f)=F(f)$ holds
for almost all $f$ with respect to the measure $\ompp$. Since the
group $\Ke_n$ is countable, for almost all $f$ the equality
$\left[g\cdot F\right](f)=F(f)$ holds for all $g\in\Ke_n$.

Consider the stabilizer $\Stab(f)\subset\Ke_n$ of a point
$f\in\Yab$. The image of the group $\Stab(f)$ under the map
$g\mapsto c(g,f)$ is the subgroup
$\Ke_a\times\Ke_b\subset\Ke\times\Ke$. Therefore, for almost all
$f\in\Yab$ we have $F(f)\in H(\pi_1)_a\otimes H(\pi_2)_b$.

For every $\zeta \in H(\pi_1)_a \otimes H(\pi_2)_b$ the set
$F^{-1}(\zeta)$ is a measurable and a $\Ke_n$-invariant $\!\mod 0$
subset. Consider an arbitrary set $\YabD$. The action of the
group $\Ke_n$ on the set $\YabD$ is clearly isomorphic to the
action of the symmetric group by permutations of factors on
$\{1,2\}^\infty$ with a Bernoulli measure. The ergodicity of
this latter action is well known. Hence the set $F^{-1}(\zeta)$
is, up to a subset of measure zero, a union of the sets $\YabD$.

Therefore, $\left(\HYab\right)_n=\Hab$.
\qed\enddemo

\proclaim{Corollary}
The representation of the group $\Ge$ in the space $H$ is
admissible. Its depth is $d=d_1+d_2$ and
$$
H_{d}=
H_{d_1,d_2} \,,
$$
$$
\dim H_{d}=\binom {2d}{2d_1} 
\dim H(\pi_1)_{d_1}
\dim H(\pi_2)_{d_2} \,.
$$
\endproclaim

\definition{Definition}
Let $\pi$ be the representation of the group $\Ge$ in the cyclic
span of  $H_{d_1,d_2}$. We call the representation $\pi$ the {\it
mixture} of representations $\pi_1$, $\pi_2$ with the weights
$p_1$, $p_2$. 
\enddefinition 

\proclaim{Proposition 10}
Suppose that the representations $\pi_1$, $\pi_2$ are irreducible
and that $\mu_i$, $i=1,2$ are their Thoma measures.
Then the  representation of the pair $\left(\Ge_d,\Ke_d\right)$ in the
cyclic span of $R(\pi)$ is a multiple of the spherical
representation with the Thoma measure
$$
\mu\left( x \right)=
p_1\mu_1\left( x/p_1 \right)+
p_2\mu_2\left( x/p_2 \right)     \,.
$$
\endproclaim

\demo{Proof}
We have to prove that
$$
\pi\left(\prod_{i}{C_i}\right)=\prod_{i}{\pi(C_i)}
$$
and
$$
\pi(C_k)=p_1^k \pi_1(C_k)+ p_2^k \pi_2(C_k)    \,.
$$
Let us prove the second equality. The proof of the first one is
similar. Denote by $z_k\in \Ge_d$ the permutation
$$
z_k=(d+1,\dots,d+k)  \,.
$$
Choose a vector $\xi\in H(\pi_1)_{d_1}\otimes H(\pi_2)_{d_2}$,
$\|\xi\|=1$. Consider the following map $F\in R(\pi)$,
$$
F(f)=\cases
\left(
\ompp\left(\Ydd \right)
\right)^{-1} \xi
,& f\in \Ydd\\
0,& \text{ otherwise }    \,.
\endcases
$$
It is clear that $\|F\|=1$. Consider two subsets
$$
{}^{(1)}\Ydd=
\left\{f\in \Ydd|
f(i)=1,|i|=d+1,\dots,d+k,
\right\}
$$
$$
{}^{(2)}\Ydd=
\left\{f\in \Ydd|
f(i)=2,|i|=d+1,\dots,d+k,
\right\}   \,.
$$
in $\Ydd$. Clearly, the two inclusions $f\in\Ydd$ and $z_k\cdot
f\in\Ydd$ occur simultaneously only  if $f\in{}^{(1)}\Ydd$ or
$f\in{}^{(2)}\Ydd$. Therefore,
$$
\gather
\pi(C_k)=\left(\pi(z_k)F,F\right)_H \\
=
\frac{1}
{\ompp\left(\Ydd \right)}
\left[
\pi_1(C_k)\,
\ompp\left({}^{(1)}\Ydd \right)+
\pi_2(C_k)\,
\ompp\left({}^{(2)}\Ydd \right)
\right]\\
=\pi_1(C_k) p_1^k+  \pi_2(C_k) p_2^k \,. \qed
\endgather
$$
\enddemo

Now consider the action of the operators $A_i$, $|i|=1,\dots,d$
in the space $R(\pi)$. Set $m_i=\dim R(\pi_i)$, $i=1,2$. Let
$\{\zeta^{(i)}_1,\dots,\zeta^{(i)}_{m_i}\}\in R(\pi_i)$ be the
eigenbases of the operators $A_k$, $|k|=1,\dots,d$. Let $D$ run
over the subsets of cardinality $2d_1$ in the set
$\{-d,\dots,d\}\setminus0$. Then the maps $\FijD$,
$i=1,\dots,m_1$, $j=1,\dots,m_2$, where
$$
\FijD(f)=
\cases
\zeta^{(1)}_i\otimes\zeta^{(2)}_j,
& f\in Y_{d_1,d_2,D}\\
0,& \text{otherwise}
\endcases
$$
form a basis in the space $R(\pi)$. This basis consists of
eigenvectors of the operators $A_k$. Indeed, let $P$ be the
orthogonal projector of the space $H$ onto the subspace $R(\pi)$.
Then by (1.5), for every $F\in R(\pi)$ and every
$k,|k|=1,\dots,n$, we have the following equality
$$
\pi(A_k)\cdot F=P \pi\left((k,d+1)\right)\cdot F   \,.
$$
Therefore
$$
\pi(A_k)\cdot\FijD(f)=
\cases
p_1 \pi_1(A_{\eta_1(k)})\zeta^{(1)}_i\otimes\zeta^{(2)},
& f\in Y_{d_1,d_2,D},k\in D\\
p_2 \zeta^{(1)}_i\otimes\pi_2(A_{\eta_2(k)})\zeta^{(2)},
& f\in Y_{d_1,d_2,D},k\notin D\\
0,& \text{otherwise}   \,.
\endcases
$$

Let $\L_1$ and $\L_2$ denote the Young distributions of
representations $\pi_1$ and $\pi_2$. Let $T_{\L_i}$, $I=1,2$,
denote the corresponding representations of the semigroups
$$
S(d_i) \ltimes \bZp^{d_i}
$$
in the spaces $R(\pi_i)$. Denote by $\L(\cdot/p_1)$ the Young
distribution which equals $\L(x/p_1)$ at the point $x$. Then the
above formulas for the action of operators $A_i$ imply that the
representation of the semigroup $S(d) \ltimes \bZp^d$ in the
space $R(\pi)$ is
$$
\Ind_{
\left( S(d_1) \ltimes \bZp^{d_1} \right)
\times
\left( S(d_2) \ltimes \bZp^{d_2} \right)
}
^{\, S(d) \ltimes \bZp^{d} }
T_{\L_1(\cdot/p_1)} \otimes T_{\L_2(\cdot/p_2)}   \,.
$$
This representation is
irreducible if and only if $\supp\Lambda_1\left(\cdot/p_1\right)\,\cap\,\supp \Lambda_2\left( \cdot /p_2 \right)=\emptyset$.
In this case the representation $\pi$ is also irreducible and
has the Young distribution
$\L$ equal to
$\Lambda\left(x\right)=\Lambda_1\left(x/p_1\right)\cup\Lambda_2\left(x/p_2\right)$.

We summarize this discussion as follows:

\proclaim{Theorem 6}
Let $\pi_1$ and $\pi_2$ be two irreducible admissible
representations of the pair
$\GKe$ with Thoma measures  $\mu_i$, $i=1,2$, and
with Young distributions $\L_i$, $i=1,2$. Let the representation
$\pi$ be the mixture of representations $\pi_1$ and $\pi_2$ with
the weights $p_1$ and $p_2$, where $p_1>0$, $p_2>0$, $p_1+p_2=1$.

The representation $\pi$ is admissible. It is a sum of
irreducible admissible representations with Thoma measure
$$
\mu\left( x \right)=
p_1\mu_1\left( x/p_1 \right)+
p_2\mu_2\left( x/p_2 \right)   \,.
$$
It is irreducible if and only if
$$
\supp \Lambda_1
\left( \cdot/p_1 \right) \cap
\supp \Lambda_2
\left( \cdot/p_2 \right)=\emptyset  \,.
$$
In this case it has the Young distribution $\L$, where
$$
\Lambda
\left( x \right)=
\Lambda_1
\left( x/p_1 \right)\cup
\Lambda_2
\left( x/p_2 \right) \,.
$$
\endproclaim

The mixture of any finite number of representations may be
defined in the same way.
Moreover, one can define a mixture of finitely many admissible
representations and countably many spherical representations
$\pi_i$, $i=1,\dots$ of the pair $\GKe$. In this case
$\bigotimes_1^\infty {H\left( \pi_i \right)}$ denotes the direct
limit of Hilbert spaces
$$
\limind \bigotimes_1^N {H\left( \pi_i \right)}
$$
with respect to the inclusions
$$
\bigotimes_1^N {H\left( \pi_i \right)}\to
\bigotimes_1^N
{H\left(\pi_i \right)}\otimes \xi_{N+1} \,,
$$
where $\xi_i\in H\left(\pi_i\right)$ is the spherical vector of the 
representation $\pi_i$. In this case Theorem 6 can be
generalized as follows.

\proclaim{Theorem 7}
Let $\pi_i$, $i=1,2,\dots$,  be irreducible admissible representations
of the pair
$\GKe$ with Thoma measures $\mu_i$, $i=1,2,\dots$, and Young
distributions $\L_i$, $i=1,2,\dots$,  such that the sum
$$
\sum_i{|\L_i|}<\infty
$$
is finite. Define the representation $\pi$ as the mixture of
representations $\pi_i$, $i=1,2,\dots$,  with the weights $p_i$,
$i=1,2,\dots$, where $p_i>0$, $\sum_i{p_i}=1$.

The representation $\pi$ is admissible. It is a sum of
irreducible admissible representations with Thoma measure
$$
\mu\left( x \right)=\sum_i
{p_i\, \mu_i\left( x/p_i \right)}   \,.
$$
It is irreducible if and only if
$$
\supp \Lambda_i \left( \cdot /p_i \right) \cap
\supp \Lambda_j \left( \cdot /p_j \right)
=\emptyset\,, \quad i\ne j\,.
$$
In this case it has the following Young
distribution:
$$
\Lambda
\left( x \right)=
\bigcup_i{
\Lambda_i
\left( x/p_i \right)} \,.
$$
\endproclaim

\subhead 3.2 Mixtures and induction \endsubhead

In this Section we show that the operation of mixing the
representations is intimately related to that of inducing of
representations.
In the group of all (not necessarily finite) bijections
$g:\bZz\to\bZz$ consider the subgroups
$$
\bKe_n=
\left\{g|g(i)=-g(-i),g(i)=i,i>n \right\}
,n=0,1,\dots \,.
$$
As usual, write $\bKe=\bKe_0$. Set
$$
\bGe=\Ge\cdot\bKe  \,.
$$
Define a topology in the group $\bGe$ in which a fundamental
neighborhood system of unity is formed by the subgroups
$\bKe_n$, $n=0,1,\dots$. In this topology the group $\Ge$ is a
dense subgroup of the group $\bGe$. The representations of the
group $\Ge$ which admit a continuous extensions onto the group $\bGe$
are exactly the admissible representations of the pair $\GKe$
\cite{32}.

A mixture of representations is an induced representation in
the following sense. The space $Y$ is a homogeneous space of the
group $\bGe$. The stabilizer of a point is isomorphic to the group
$\bGe\times\bGe$. The cocycle $c(\cdot,\cdot)$ is exactly the
usual cocycle on the homogeneous space. We have the following

\proclaim{Proposition 11}
The measures $\ompp$ are, up to a factor, exactly  the measures on
$Y$ which are
\roster
\item invariant with respect to the action of the group $\Ge$;
\item finite on all the sets $Y_n$, $n=0,1,\dots$;
\item extreme in the class of measures with the properties
(1)--(2).
\endroster
\endproclaim

\demo{Proof}
Consider the set $Y_k$. For every function $f\in Y_k$ there is an
element $g\in\Ge(k)$ such that $g\cdot f\in Y_0$. In other words,
$Y_k\subset\bigcup_{g\in\Ge(k)}{g\cdot{Y}_0}$. This means that
a $\Ge$-invariant measure on $Y$ is determined by its restriction to
$Y_0$.

Let $\nu$ be a measure satisfying the conditions (1)--(3). If
$\nu(Y_0)=0$, then also $\nu(Y)=0$. Therefore we can assume that
$\nu(Y_0)=1$. The set $Y_0$ is embedded in
$X_0\cong\{1,2\}^\infty$. The group $\Ke$ acts on $X_0$ by
permutations of the factors. It follows from the de Finetti
theorem \cite{53, p.~256} that the measure $\nu$ on $X_0$ has the
form $\ompp$ for some $p_1$, $p_2\geq0$, $p_1+p_2=1$. It is also
clear that $\nu(Y_0)\neq0$ if $p_1>0$ and  $p_2>0$.
\qed\enddemo

Let us now consider the $\Ge$-invariant measures supported
on $\bY\setminus Y$, where 
$\bY$ is the closure of the set $Y$ in the space $X$: 
$$
\bY=\left\{f,
|f^{-1}(1)\cap \{-N,\dots,N\}|
\text{ is even for almost all }N \right\}    \,.
$$
The set $\bY$ coincides with the orbit of the set $X_0$ under the
action of the group $\Ge$. The set $\bY\setminus Y$ consists of a
countable number of $\Ge$-orbits $Z^{(i)}_{2k}$,
$i=1,2,\,k=0,1,\dots$
$$
Z^{(i)}_{r}=\left\{f\in X,|f^{-1}(i)|=r \right\}  \,.
$$
It suffices to consider the sets
$Z^{(1)}_{2k}$. The set $Z^{(1)}_{2k}$ is a single orbit
of a countable group $\Ge$. It supports a unique, up to a
factor, $\Ge$-invariant measure $\nu_{2k}$ which  is just
the counting measure. The stabilizer of a point in $Z^{(1)}_{2k}$ is
isomorphic to the group $\Ge(k)\times\Ge_k\cong \Ge(k)\times\Ge$. 

Let $\rho$ be an irreducible representation of the group
$\Ge(d_1)$, and $\pi_2$ be an irreducible admissible
representation of the group $\Ge_{d_2}\cong \Ge$. Let $\mu_2$,
$\L_2$ be the Thoma measure and the Young distribution
corresponding to the representation $\pi_2$. The representation
$\pi$,
$$
\pi=
\Ind_{\,\Ge(d_1)\times\Ge_{d_2}}
^{\,\Ge} \rho \otimes \pi_2\,,
$$
can be realized  in the space of maps $F:X\to H(\rho)\otimes H(\pi_2)$
which are square summable with respect to the measure $\nu_{2d_1}$.

Assume that $\L_2(0)=\emptyset$. Then, the same argument we used for
mixtures yields that the
representation $\pi$ is admissible and irreducible. Its Thoma
measure is
$$
\mu=\mu_2  \,,
$$
and the Young distribution is
$$
\L(x)=
\cases
\L_2(x),& x\neq 0 \\
\rho,& x=0 \,.
\endcases
$$

This conventional induction of representations may be considered
as a limit case of mixtures. Recall the definition of the topology
in the space of unitary representations of a discrete group $G$
\cite{23}. Let $T_0$ be a unitary representation. Given a finite
subset $M\in G$ and an array of vectors
$\xi_1,\dots,\xi_k\in{H}(T_0)$, we denote by
$U(T_0,M,\xi_1,\dots,\xi_k;\varepsilon)$ the set of unitary representations
$T$ of the group $G$, for which the corresponding space contains
the vectors $\zeta_1,\dots,\zeta_n$ such that
$$
|(T(g)\zeta_i,\zeta_j)-(T_0(g)\xi_i,\xi_j)|<\varepsilon\,,\quad g\in M,
i,j=1,\dots,n \,.
$$
The sets $U(T_0,M,\xi_1,\dots,\xi_k;\varepsilon)$ form a neighborhood base
of the representation $T_0$. If the representation $T_0$ is
irreducible then in order to check the convergence $T_n\to T_0$
it suffices to check that a certain matrix element
$(T_0(\cdot)\xi,\xi),\xi\in H(T_0)$ can be approximated by the
matrix elements of representations $T_n$.

Let $\pi_1$ be an irreducible admissible representation of
the group $\Ge$ with the Young distribution $\L_1$ such that
$\supp \L_1$ is a one point set $\{y\}$, $y\in [-1,1]$ and
$\L_1(y)=\rho$. Denote by $\pi(p_1,p_2)$ the mixture of
representations $\pi_1$, $\pi_2$ with the weights $p_1$, $p_2$.

\proclaim{Proposition 12}
$$
\pi(p_1,p_2) \to \pi \text{  as  } p_1 \to 0    \,.
$$
\endproclaim

\demo{Proof}
It follows from Theorem 6 that, as $p_1\to0$, the representations
of the semigroup $\Gamma(d)$ in the spaces $R(\pi(p_1,p_2))$
converge to the representation of this semigroup in $R(\pi)$.
Hence, the matrix coefficients of representations $\pi(p_1,p_2)$
corresponding to vectors in subspaces $R(\pi(p_1,p_2))$ converge
to the matrix coefficients of the representation $\pi$. Since the
representation $\pi$ is irreducible, the proposition follows.
\qed\enddemo

The indecomposable invariant measures on $X$ supported by the set
$X\setminus\bY$ have the following meaning. The restriction of a
measure $\ompp$ onto $X\setminus\bY$ corresponds to the
representations of the group $\Ge$ which are the mixtures of two
representations of the group $\Go$ with the weights $p_1$, $p_2$.
The measures $\nu_{2k+1}$ correspond to representations of the
group $\Ge$ which are induced from a subgroup isomorphic to
$\Go(k)\times\Go$.

\subhead 3.3 Elementary representations \endsubhead

Let $\pi$ be an irreducible admissible representation of $\Ge$ or
$\Go$ with Thoma measure $\mu$ and Young distribution $\L$.
Call this representation {\it elementary} if
$$
\supp\mu=\supp\L=\{y\}
$$
for some point $y\in [-1,1]$. It follows from the classification
of irreducible representations and the results of the two
previous sections that in order to construct all irreducible
admissible representation it remains to give a construction
of elementary ones.

We shall
briefly describe the realization of elementary representations of
the groups $\Ge$ and $\Go$ obtained in \cite{32}. One has to
distinguish between three cases: $y=0$, $y>0$ and $y<0$.

Suppose $y=0$. Then the corresponding elementary representation
is
$$
\Ind_{\,G(d)\times K_d}^{\,G} \L(0)\otimes 1   \,.
$$

Suppose $y>0$. By Theorem 4, $y$ has to be of the form
$y=1/n$, $n\in\bN$. Consider the space $\bC^n$ with 
the standard basis $\{e_1,\dots,e_n\}$ and
set, by definition,
$$
\xi=n^{-1/2}\sum_{i=1}^n{e_i\otimes e_i} \in \bC^n\otimes\bC^n \,.
$$
Consider the direct limit
$$
H=\limind \bigotimes_1^{2N}{\bC^n}
$$
of Hilbert spaces with respect to inclusions
$$
\bigotimes_1^{2N}{\bC^n}\to
\bigotimes_1^{2N}{\bC^n}\otimes\xi    \,.
$$
Using the bijection $\bZz\to\bN$,
$$
i\mapsto
\cases 2i-1,&i>0 \\
-2i,&i<0\,,
\endcases
$$
we can define the action of the group $\Ge$ in the space $H$. (In
case of $\Go$ one should consider
$\limind\bigotimes_0^{2N}{\bC^n}$.)

The group of orthogonal matrices $O(n)$ preserves the vector
$\xi$, hence its action on the space $H$ is well defined. As
explained in \cite{54}, the irreducible 
representations of the group $O(n)$
are labeled by Young diagrams $\lambda$ such that
$$
(\lambda)'_1+(\lambda)'_2 \leq n   \,.
$$
One has the following:

\proclaim{Theorem \rm (Olshanski, \cite{32})}
\roster
\item The representation of the group $\Ge$ in the space $H$ is
admissible.
\item The representations of the groups $\Ge$ and $O(n)$ generate
the commutant of each other.
\item The space $H$, as an $\Ge\times O(n)$-module, decomposes into
the following  direct sum
$$
H=\bigoplus_
{\lambda,\,\,(\lambda)'_1+(\lambda)'_2 \leq n}
{\pi_{1/n,\lambda}\otimes T_{\lambda}}   \,,
$$
where $T_{\lambda}$ is the representation of the group $O(n)$
corresponding to a diagram $\lambda$, and $\pi_{1/n,\lambda}$ is
the irreducible admissible representation of the group $\Ge$, such
that $\supp\mu\!=\!\supp\L\!=\!\{1/n\}$, $\L(1/n)\!=\!\lambda$.
\endroster
\endproclaim

There is a little inaccuracy in the paper \cite{32} in case of
$y<0$, indicated by G.~Olshanski. The correct construction of
representations is as follows. By Theorem 4, in case of
$y<0$ we are forced to take $y=-1/2n$, $n\in\bN$. Consider a
basis $\{e_1,\dots,e_{2n}\}$ in the space $\bC^{2n}$. Consider
the vector
$$
\xi=(2n)^{-1/2}\sum_{i=1}^{n}
{(e_i\otimes e_{i+n}-e_{n+i}\otimes e_i)}
$$
in the space $\bC^{2n}\otimes\bC^{2n}$ and let
$$
H=\limind \bigotimes_1^{2N}{\bC^{2n}}
$$
be the direct limit of Hilbert spaces with respect to inclusions
$$
\bigotimes_1^{2N}{\bC^{2n}}\to
\bigotimes_1^{2N}{\bC^{2n}}\otimes\xi     \,.
$$
Define the representation of the group $\Ge$ in the space $H$ as
the tensor product of the representation by permutations of tensor
factors, and the one dimensional representation $\sgn$.

The group $Sp(n)$ of symplectic matrices preserves the vector
$\xi$, hence its action on the space $H$ is well defined.
Again, as explained in \cite{54}, the representations of
the group $Sp(n)$ are labeled by Young diagrams $\lambda$, such
that
$$
\lambda_1\leq n   \,.
$$
One has the following

\proclaim{Theorem \rm (Olshanski)}
\roster
\item The representation of the group $\Ge$ in the space $H$ is
admissible.
\item The representations of groups $\Ge$ and $Sp(n)$ generate
the  commutant of each other.
\item The space $H$, as an $\Ge\times Sp(n)$-module, decomposes
into the following direct sum
$$
H=\bigoplus_
{\lambda,\,\,\lambda_1\leq n}
{\pi_{-1/2n,\lambda}\otimes T_{\lambda}}   \,,
$$
where $T_{\lambda}$ is the representation of the group $Sp(n)$
corresponding to a diagram $\lambda$, and $\pi_{-1/2n,\lambda}$
is an irreducible admissible representation of the group $\Ge$,
such that $\supp\mu\!=\!\supp\L\!=\!\{-1/2n\}$,
$\L(1/2n)\!=\!\lambda$.
\endroster
\endproclaim

\subhead 3.4 Mixtures in the case of $\Gd$ \endsubhead

The definition of a mixture of representations can be easily
extended to the case of the pair $\GKd$. In fact, it is natural
to consider the following  ``unbalanced'' groups
$$
\Gd_{m_1,m_2}=
\left\{g\in \Gd | g(i)=i, -m_2 \leq i \leq m_1
\right\}    \,,
$$
where $m_1,m_2$ are integers, and define
mixtures of representations of these groups. Since the inclusion
$\Kd_n\subset\Gd_{m_1,m_2}$ is valid, given $m_1,m_2$, for all
$n$ but finitely many, the definition of admissible
representation also works for the groups $\Gd_{m_1,m_2}$. All the
theory of admissible representations can be transferred
word-for-word to this ``unbalanced'' case. In particular, such
representations are labeled by a Thoma measure and a pair of
Young distributions $\L,\Mu$, for which it is now possible that
$|\L|\neq|\Mu|$.

The construction of elementary representations can be taken from
\cite{32}.

Remark that in the language of $H$-series from Section 2.4 the mixture
of spherical representations corresponds to the product of
$H$-series. 

\head 4. Concluding remarks
\endhead

This paper is the English version of the author's PhD thesis
(1995, Moscow State University). I did not try to update 
anything in it. For some related recent results the reader is
referred to \cite{57,58} and references therein.

I am very grateful to my advisor A.~A.~Kirillov, A.~M.~Vershik,
R.~S.~Ismagilov, S.~V.~Kerov, and Yu.~A.~Neretin for their
constant interest, encouragement, and help. My very special
thanks are due to G.~Olshanski, who not only laid the
foundations of the whole subject, thus making this paper
possible in the first place, but also was of absolutely
indispensable help to me from the very beginning of my work 
to the proofreading stage of the present paper. 


\enddocument
\end